\documentclass[preprint,12pt]{elsarticle}
\usepackage{graphicx} 
\usepackage{subcaption}
\usepackage{nomencl}
\usepackage{etoolbox}
\usepackage{amssymb}
\usepackage{hyperref}
\usepackage[utf8]{inputenc}
\usepackage[T1]{fontenc}
\usepackage{float}
\usepackage{color}
\usepackage{lineno}
\usepackage{amsmath}
\usepackage[linesnumbered]{algorithm2e}
\usepackage{caption}
\usepackage{booktabs}
\usepackage{array}
\usepackage{pgfplots}
\usepackage[euler-digits]{eulervm}
\usepackage[eulergreek]{sansmath}
\usepackage{amsmath} 
\usepackage[normalem]{ulem}

\usepackage{listings}

\lstdefinestyle{custompython}{
  belowcaptionskip=1\baselineskip,
  breaklines=true,
  frame=L,
  xleftmargin=\parindent,
  language=python,
  showstringspaces=false,
  basicstyle=\footnotesize\ttfamily,
  keywordstyle=\bfseries\color{green!40!black},
  commentstyle=\itshape\color{purple!40!black},
  identifierstyle=\color{blue},
  stringstyle=\color{orange},
}
\usetikzlibrary{positioning}
\pgfplotsset{compat=1.18}
\newlength\algowd

\definecolor{basquered}{RGB}{213,43,30}
\newcommand{\revision}[1]{\textcolor{black}{#1}}

\begin{document}

\begin{frontmatter}

\title{
Augmenting MRI scan data with real-time predictions of glioblastoma brain tumor evolution using  \revision{faster} exponential time integrators}

\author{Magdalena Pabisz$^1$, Judit Muñoz-Matute$^{2,3}$, Maciej Paszy\'nski$^1$}

\address{$^1$ AGH University of Krakow, Poland\\  $^2$ The Basque Center for Applied Mathematics (BCAM), Bilbao, Spain \\ $^3$ The Oden Institute for Computational Engineering and Sciences, \\ The University of Texas at Austin, USA}

\begin{abstract}
We present a MATLAB code for exponential integrators method simulating the glioblastoma tumor growth. \revision{It employs the Fisher-Kolmogorov diffusion-reaction tumor brain model with logistic growth. The input is the MRI scans of the human head and the initial tumor location. The simulation uses the finite difference formulation in space and the ultra-fast exponential integrators method in time. The output from the code is the input data for ParaView visualization.
While there are many brain tumor simulation codes, our method's novelty lies in its implementation using exponential integrators.
We propose a new algorithm for the fast computation of exponential integrators. Regarding execution time on a laptop with Win10, using MATLAB, with 11th Gen Intel(R) Core(TM) i5-11500H, 2.92 GHz, and 32 GB of RAM, the algorithm outperforms the state-of-the-art routines from \cite{al2011computing}.
We also compare our method with an implicit, unconditionally stable Crank-Nicolson time integration scheme based on the finite difference method. We show that our method is two orders of magnitude faster than the Crank-Nicolson method with finite difference discretization in space on a laptop equipped with MATLAB.
The brain tumor two-year future prediction using $128\times 128 \times 128$ computational grid and 100-time steps, built over the MRI scans of the human head, takes less than 10 minutes on the laptop. }
\end{abstract}


\end{frontmatter}

\section{Introduction}

\revision{Tumor growth numerical models are essential for oncology research and clinical applications. They provide a quantitative framework for understanding how tumors grow, spread, and respond to treatment.
Fast tumor models can provide real-time tumor progression or regression predictions in response to various therapies. Accurate predictions are vital to avoid over- or under-treatment and optimizing patient outcomes. In personalized medicine, fast and precise tumor models can simulate the effects of different therapies on a specific patient's tumor, allowing clinicians to choose the most effective treatment.
Fast models can be used to optimize the dosing, timing, and delivery of different therapies. The therapy optimization requires repeated simulations to find the best possible treatment plan, necessitating both speed and accuracy.}

\revision{Fast and accurate tumor growth numerical models are crucial for advancing cancer research and treatment. Their speed enables real-time decision-making, optimization of treatment plans, and rapid hypothesis testing, while their accuracy ensures that predictions are reliable and clinically relevant. By balancing these requirements, these models can significantly impact personalized medicine, drug development, and our understanding of cancer biology.}

Glioblastoma is a malignant brain tumor with a high mortality rate \cite{rate}.
This tumor is highly aggressive, and it generates microvascular proliferation that is not visible on MRI scans \cite{a5c33ecf60d74369a20b8fe472f6ceea}. Thus, computer-based simulation predictions of the evolution of the brain tumor are essential in planning treatment and surgery. But highly efficient simulators of the tumor growth take several minutes to compute a single time step \cite{Lipkov2018PersonalizedRD,LOS20171257,LOS20191,doi:10.1177/1094342018816772}.

\revision{An overview of the modern computational methods for tumor growth simulations are described in \cite{GPU3}.
While the two-dimensional tumor growth simulations \cite{JOCS1,JOCS2} are relatively fast and do not require unique acceleration methods, GPUs have been employed for speeding up 3D tumor growth simulations for over a decade \cite{GPU4}.
The current state-of-the-art in HPC tumor simulations can be summarized in the following way.
The finite difference 3D tumor growth solver from \cite{doi:10.1177/1094342018816772} using a system of advection-diffusion-reaction equations, can solve one time step on $250\times 250\times 250$ stencil points using one GPU within 100s, or one time step on $750\times 750\times 750$ stencil points using 8 GPU within 300s on Prometheus cluster.
The GPU simulations \cite{GPU1} of a simplified (using one PDE) model of 3D tumor growth on the computational mesh of $240\times 120\times 120$ points take 39s per one-time step on GPU versus 1641s (half an hour) on CPU.
The higher-order finite element method solver takes 2122s on CPU or 12s on two GPUs to compute a one-time step on $256\times 256\times 256$ mesh using IGA-FEM solver \cite{GPU2}.
There are also hybrid simulators, e.g., coupling the finite elements with chemical interactions \cite{Parallel2}, computing 9,000-time steps within six days, incorporating 75,000 interacting fibers into finite element method solver, or cellular automata massive parallel tumor growth simulators \cite{Parallel1} executed on 2040 cores of Intel Xeon Phi processors or 30 computing nodes of Stampede cluster, reaching 55 speedup.}
The Physics Informed Neural Networks originally proposed by Karniadakis et al. \cite{raissi2019physics,chen2020physics,cai2021physics} instantiated for the brain tumor simulation in \cite{Zhu} can predict the future behavior of the glioblastoma tumor evolution within one hour for the patient-specific case.

\revision{In this paper, we will show how to employ the exponential integrators method to compute one time step on $128\times 128\times 128$ stencil points within 10 seconds on a laptop (two orders of magnitude faster than alternative state-of-the-art finite difference or finite element methods).}

We present a computational tool to perform ultra-fast simulations of two and three-dimensional progression of glioblastoma brain tumor models. We consider the Fisher-Kolmogorov diffusion-reaction equation with logistic growth \cite{Logistic} describing the evolution of glioblastoma tumor in a non-homogenous environment of the human head. Our simulations employ the MRI scan data consisting of 29 bitmaps, with pixel intensity varying from 0 to 255.
\revision{Based on the pixel intensity from the bitmaps, we decide what kind of material data we have there. We introduce skull, skin, white matter, gray matter, and air based on the predefined thresholds of the pixel intensities.}

We combine two methods to discretize the transient semilinear model. First, we semidiscrete the space variable employing finite differences.
\revision{Then, we employ the exponential time integrators \cite{hochbruck2010exponential,munoz2023multistage,CROCI2023101966} designed to solve semilinear problems. They are suitable for long-time simulations and unconditionally stable.}
For simplicity, we focus on the exponential Euler method, which is first-order, but extending to higher-order exponential-based methods in time is straightforward.
\revision{We propose a new algorithm for fast computations of the exponential integrators. Our MATLAB implementation of the exponential integrators method is faster than the state-of-the-art MATLAB library \cite{al2011computing}
We also compare our exponential integrators simulator to the implicit dynamics simulation using the unconditionally stable Crank-Nicolson method and finite difference discretization in space. Using MATLAB implementation, our method is up to two orders of magnitude faster than Crank-Nicolson method on a laptop equipped with Win10, using MATLAB, with 11th Gen Intel(R) Core(TM) i5-11500H, 2.92 GHz, and 32 GB of RAM}.

Then, we post-processed the results and used ParaView for visualization.
We test this strategy on 2D and 3D non-stationary glioblastoma models employing MRI data for the diffusion coefficient. The numerical results show that we can perform $100$ iterations over
\revision{$128\times128\times128$ computational stencil points in less than 10 minutes on the laptop}. Thus, this simulator can be employed "on the fly".

\revision{While there are many brain tumor simulation codes \cite{ReviewTumors}, the novelty of our method lies in a new implementation using a new algorithm for computing the exponential integrators method.}
From the authors understanding, there is only one work in the literature \cite{maddalena2020existence} where the authors employ exponential time integrators to simulate cancer models. Our article aims to underscore the potential of this tool in the simulation of glioblastomas in the human brain. \revision{In particular, we show that our method is up to two orders of magnitude faster than alternative state-of-the-art methods based on the Crank-Nicolson scheme and finite difference discretization and several times faster than the existing state-of-the-art exponential integrator library \cite{al2011computing}.}

Our future work will include developing the parallel version of the simulator and targeting GPGPU cards. This will allow us to perform simulations on high-resolution data from the brain. We also plan to incorporate drug delivery into our model \cite{McDaniel2013}. With the chemotherapy in our simulator, we plan to attack the data assimilation problem using the supermodeling technique \cite{PASZYNSKI2022214,SIWIK2022114308}.

\section{Glioblastoma tumor model}

\begin{table}[]
    \centering
    \begin{tabular}{|c|c|c|}
         \hline
      \revision{Variable} &  \revision{Meaning} &  \revision{Value} \\
\hline
      \revision{$u(x,y,z;t)$} &  \revision{tumor cell density} & \revision{[0,1]} \\
      \revision{$\hat{D}$} & \revision{"global" diffusion estimate} & \revision{$p_wD_w+p_gD_g$} \\
      \revision{$p_w$} & \revision{proportion of white matter} & \revision{70\%(2D) 80\%(3D)} \\
      \revision{$p_g$} & \revision{proportion of gray matter}  & \revision{30\%(2D) 20\%(3D)} \\
      \revision{$D_w$} & \revision{isotropic diffusion in white matter} & \revision{0.13 [mm$^2$days$^{-1}$]}\\
      \revision{$D_g$} & \revision{isotropic diffusion in gray matter}  & \revision{0.013 [mm$^2$days$^{-1}$]}\\
      \revision{$D(x,y,z)$} & \revision{"local" diffusion at point $(x,y,z)$} & \revision{$D_w$ or $D_g$ or 0}\\
      \revision{$\rho$} & \revision{proliferation rate of tumor cells}  & \revision{0.025 [days$^{-1}$]}\\
      \revision{$t_0$} & \revision{initial time moment} & \revision{150 days} \\
	\revision{$T$} & \revision{final time moment} & \revision{3500 days} \\
      \revision{$\tau$} & \revision{time step size} & \revision{$\frac{3500-150}{100}=33$ [days]} \\
   	\revision{$K$} & \revision{number of time steps} & \revision{100} \\
         \hline
    \end{tabular}
    \caption{\revision{Variables of the glioblastoma tumor model.}}
    \label{tab:times}
\end{table}

Let $\Omega\subset {\mathbb{R}}^3$ an open set and $I=(0,T)$ with $T>0$, we consider the following Fisher-Kolmogorov diffusion-reaction equation with logistic growth \cite{Logistic} describing the brain tumor dynamics.
\begin{equation}\label{PDE}
\displaystyle{ \left\{
\begin{aligned}
\frac{\partial u}{\partial t}=\underbrace{\nabla \cdot(D(\mathbf{x})\nabla u)}_{\textrm{Tumor cell diffusion}}&+\underbrace{\rho u(1-u)}_{\textrm{Tumor cell proliferation}},&\mbox{in}&\;\Omega \times I,\\
\nabla u\cdot n &= 0,&\mbox{on}&\;\partial \Omega \times I,\\
u(\mathbf{x},0)&=u_{0},&\mbox{on}&\;\Omega \times \{0\}.\\
\end{aligned}
\right.} 
\end{equation}
Here, $u(\mathbf{x};t)$ represents the tumor cell density scalar field, and $D(\mathbf{x})$ represents the diffusion coefficient, \revision{one for each kind of tissue},
describing how the tumor cells expand in a tissue, similar to the diffusion phenomena.
The $\rho$ coefficient describes the proliferation rate of the tumor cells (how fast they multiply).

\revision{Our model is heterogeneous.
The diffusion $D_{i,j,k} = D(x_{i,j,k})$ changes with the points of the computational stencil. For the white matter, it is equal to 0.13; for the gray matter, it is equal to 0.013; and for the air (outside the head), it is equal to 0.
There is no tumor diffusion in the air outside the human head.}

\revision{On the other hand, the "global" diffusion coefficient $\hat{D}$ represents the global diffusion resulting from a proportion between white and grey matter, and it is employed for estimating the theoretical tumor growth velocity.}
Modeling the tumor cell dynamics by this kind of PDE results in a progression of the tumor "traveling wave" with the velocity $2\sqrt{\hat{D}\rho}$, according to \cite{https://doi.org/10.1046/j.1365-2184.2000.00177.x}.
\revision{The global diffusion coefficient usually is patient-specific}. Following \cite{c25}; it is expressed as a linear combination of gray and white matter coefficients, namely $\hat{D}=p_wD_w+p_gD_g$, where $p_w,p_g$ corresponds to the proportions of the white and gray matter, obtained from the MRI scans. Here, we considered the case of isotropic diffusion in both the white and the grey matter \cite{hogea2008image}. Following \cite{Menze} and \cite{Zhu}, we assume that $D_w=10D_g$, and we select $D_w=0.13$ [mm$^2$days$^{-1}$], $D_g=0.013$ [mm$^2$days$^{-1}$]. Moreover, following \cite{c25} we set $\rho=0.025$ [days]$^{-1}$.

\subsection{\revision{Spatial discretization with finite difference method}}

\revision{
We consider a regular computational grid with points 
\begin{eqnarray}
\{ x_{i,j,k} = \left((i-1)h,(j-1)h,(k-1)h\right), i=1,...,N_x; j=1,...,N_y,k=1,...,N_z\} 
\end{eqnarray}
where $N_x,N_y,N_z$ denotes the number of grid points in three spatial dimensions, and $h$ is the uniform grid diameter in all spatial directions.
Spatial discretization involves the finite difference method
\begin{eqnarray}
\begin{aligned}
\frac{\partial u^t_{i,j,k}}{\partial t} =& D_{i,j,k}\frac{u^t_{i+1,j,k}-2u^t_{i,j,k}+u^t_{i-1,j,k}}{h^2} \\
+& \frac{D_{i+1,j,k}\left(u^t_{i+1,j,k}-u^t_{i,j,k}\right) - D_{i,j,k}\left(u^t_{i,j,k}-u^t_{i-1,j,k}\right)}{h^2} \\
+&
D_{i,j,k}\frac{u^t_{i,j+1,k}-2u^t_{i,j,k}+u^t_{i,j-1,k}}{h^2} \\
+& \frac{D_{i,j+1,k}\left(u^t_{i,j+1,k}-u^t_{i,j,k}\right) - D_{i,j,k}\left(u^t_{i,j,k}-u^t_{i,j-1,k}\right)}{h^2} \\
  + &
D_{i,j,k}\frac{u^t_{i,j,k+1}-2u^t_{i,j,k}+u^t_{i,j,k-1}}{h^2} \\
+& \frac{D_{i,j,k+1}\left(u^t_{i,j,k+1}-u^t_{i,j,k}\right) - D_{i,j,k}\left(u^t_{i,j,k}-u^t_{i,j,k-1}\right)}{h^2} \\
+ & \rho u^t_{i,j,k}(1-u^t_{i,j,k}),
\end{aligned}
\end{eqnarray}
where $h$ is the stencil dimension, $ u^t_{i,j,k}=u(x_{i,j,k};t)$, and $ D_{i,j,k}=D(x_{i,j,k})$.
The equation in the matrix form reads
\begin{eqnarray}
\begin{aligned}
{\bf \frac{\partial u}{\partial t}}^t =& {\bf A} {\bf u}^t + {\bf F}{\bf u}^t, \end{aligned}
\end{eqnarray}
where
\begin{eqnarray}\label{Umatrix}
\{ {\bf \frac{\partial u}{\partial t}}^t \}_{i,j,k} = 
\frac{\partial u^t_{i,j,k}}{\partial t}, \qquad
\{ {\bf u }^t \}_{i,j,k} = u^t_{i,j,k}.
\end{eqnarray}
Using these integer coordinates, we define
\begin{eqnarray}\label{Amatrix}
h^2\{{\bf A}\}_{i,j,k;l,m,n}=\displaystyle{\left\{
\begin{aligned}
-9D_{i,j,k}-D_{i+1,j,k}-D_{i,j+1,k}-D_{i,j,k+1},  \quad (i,j,k)==(l,m,n),\\
D_{i,j,k},  \quad (l,m,n)\in \{(i+1,j,k),(i,j+1,k),(i,j,k+1)\},\\
2D_{i,j,k},  \quad (l,m,n)\in \{(i-1,j,k),(i,j-1,k),(i,j,k-1)\},\\
D_{i+1,j,k},  \quad (l,m,n)==(i+1,j,k),\\
D_{i,j+1,k},  \quad (l,m,n)==(i,j+1,k),\\
D_{i,j,k+1},  \quad (l,m,n)==(i,j,k+1),\\
0,  \quad \textrm{otherwise},
\end{aligned}
\right.} 
\end{eqnarray}
We employ $\{1,...,N_x\}\times \{1,...,N_y\} \times \{1,...,N_z\} \ni(i,j,k) \rightarrow global(i,j,k)=i+(j-1)N_y+(j-1)(k-1)N_yN_z$, the mapping from the integer coordinates into the global rows/columns numbering. }
\revision{We also notice that
\begin{eqnarray}
\begin{aligned}
{\bf F}(u^{t})=&\rho u^{t}(1-u^{t})\\=&\rho u^{t}-\rho (u^{t})^2\\=&\rho u^{t}-2\rho u^{t}u^{t-1}+\rho (u^{t-1})^2-\rho (u^{t-1}-u^{t})^2\\=&\rho (1-2u^{t-1})u^{t}+\rho (u^{t-1})^2-\rho (u^{t-1}-u^{t})^2, 
\end{aligned} \label{eq:F}
\end{eqnarray}
so we include the linear term $\rho \left(1-2u^{t}\right)u^{t+1}$ into the matrix ${\bf A}$, as well as define 
${\bf F}(u^{t})=\rho (u^{t-1})^2-\rho (u^{t-1}-u^{t})^2$.}

\revision{So the  ${\bf F}$ operator is given by 
\begin{eqnarray}\label{Umatrix}
\{ {\bf F}(u^t) \}_{i,j,k} = \rho (u^{t-1}_{i,j,k})^2-\rho (u^t_{i,j,k}-u^{t-1}_{i,j,k})^2.
\label{eq:newF}
\end{eqnarray}
In (\ref{eq:newF}) the contribution of term $\rho (u^{t-1})^2$ can be calculated exactly.
The definition of the right-hand side (\ref{eq:newF}) allows to reduce the discretization error in comparison to the original formulation (\ref{eq:F}).
With the original defintion (\ref{eq:F}) for $u^{t}=u^{t-1}$ we get ${\bf F}(u^{t-1}) = \rho u^{t-1}(1-u^{t-1})=\rho(u^{t-1}-(u^{t-1})^2)$. With the new definition (\ref{eq:newF})  we get
${\bf F}(u^{t-1}) = \rho (u^{t-1})^2-\rho (u^{t-1}-u^{t-1})^2 = \rho (u^{t-1})^2$.}


\subsection{\revision{Unconditionally stable exponential integrator time marching method}}

\revision{For the system of equation
\begin{equation}
\frac{d{\bf u}}{dt}={\bf A}{\bf u} +{\bf f}\left({\bf u}(t)\right)
\end{equation}
we introduce the solution as an integral equation,
an update from time moment $t_n$ into the time moment $t_{n+1}$, for time interval $I=(t_n,t_{n+1})$, where $t_{n+1}=t_n+\tau$.
In particular, we start from formula (2.2) from \cite{al2011computing}, with $u_0$ replaced by $u_n$ (previous time step solution), $t_0$ replaced by $t_n$, $t$ replaced by $t_{n+1}$. We employ $s$ as the integration variable, and $\tau=t_{n+1}-t_n$ in our paper. The equation (2.2) from \cite{al2011computing} transforms into
\begin{equation}
{\bf u}^{n+1} = \exp\left(\tau{\bf A}\right){\bf u}^{n}
+  \tau  
\int_0^1 \exp\left(\tau {\bf A}(1-s)\right){\bf f}\left({\bf u}(t_n+\tau s)\right)ds.
\end{equation}
We will project the source, and we will approximate
${\bf f}(u_n)\approx {\bf f}\left({\bf u}(t_n+\tau s)\right)$ (like in the Euler method). 
After the approximation we can take ${\bf f}(u_n)$ out, since it does not depend on $s$,
\begin{equation}
\begin{aligned}
{\bf u}^{n+1} =& \exp\left(\tau{\bf A}\right){\bf u}^{n}
+  \tau {\bf f}\left({\bf u}_n\right)  
\int_0^1 \exp\left(\tau {\bf A}(1-s)\right)ds= \\
& \exp\left(\tau{\bf A}\right){\bf u}^{n}
+  \tau {\bf f}\left({\bf u}_n\right)  
\left( -\tau {\bf A}^{-1} \left({\bf I}-\exp\left(\tau {\bf A}\right)\right)\right),
\end{aligned}
\end{equation}
we employ the property $\exp({\bf 0})={\bf I}$, with ${\bf 0}$ and ${\bf I}$ being  zero and identity matrices, repectively.
By introducing
\begin{eqnarray}
\begin{aligned}
\phi_0\left(\tau {\bf A}\right) & =\exp\left(\tau  {\bf A}\right), \\
\phi_1\left(\tau {\bf A}\right) & = \tau^{-1} {\bf A}^{-1} \left( \exp\left( \tau{\bf A}\right)-{\bf I}\right), 
\end{aligned}
\end{eqnarray}
we have
\begin{eqnarray}
\begin{aligned}
\phi_1\left(\tau {\bf A}\right) & =\tau^{-1} {\bf A}^{-1} \left( \phi_0\left(\tau {\bf A}\right)-{\bf I}\right), \\
\phi_0\left(\tau {\bf A}\right) & =\tau{\bf A}  \phi_1\left( \tau{\bf A}\right)+{\bf I}, 
\end{aligned}
\end{eqnarray}
and
\begin{equation}
{\bf u}^{n+1} = \phi_0\left(\tau {\bf A}\right){\bf u}^{n}
+  \tau {\bf f}\left({\bf u}_n\right)  
\phi_1\left(\tau {\bf A}\right),
\end{equation}
\begin{equation}
{\bf u}^{n+1} = \left(\tau{\bf A}  \phi_1\left( \tau{\bf A}\right)+{\bf I}\right){\bf u}^{n}
+  \tau {\bf f}\left({\bf u}_n\right)  
\phi_1\left(\tau {\bf A}\right).
\end{equation}
\begin{equation}
{\bf u}^{n+1} = {\bf u}^{n} + \tau \phi_1\left( \tau{\bf A} \right)\left( {\bf A}  {\bf u}^{n}
+   {\bf f}\left({\bf u}_n\right)  \right)
,
\end{equation}
We will need to compute the exponential
\begin{eqnarray}
\begin{aligned}
\phi_1\left(\tau{\bf A}\right) & =\sum_{k=1}^{\infty}\frac{\left(\tau{\bf A}\right)^{k+1}}{(k+1)!}.
\end{aligned}
\end{eqnarray}
Actually, instead of computing the exponential itself, we will compute the action of the exponential $\tau \phi_1\left(\tau{\bf A}\right) {\bf b}$
with ${\bf A} \in {\cal M}^{n \times n}$, ${\bf b} = \left( {\bf A}  {\bf u}^{n}
+   {\bf f}\left({\bf u}_n\right)\right) \in {\cal M}^{n \times 1}$. 
We will compute this action using the method summarized in \cite{al2011computing}.
Following \cite{al2011computing}, we construct a matrix ${\bf B}=\begin{bmatrix}{\bf A} & {\omega} \\ {\bf 0} & {\bf J} \end{bmatrix} \in {\cal R}^{(n+2) \times (n+2)}$, where ${\bf J}=\begin{bmatrix} 0 & 1 \\ 0 & 0 \end{bmatrix} \in {\cal R}^{2 \times 2}$, and ${\omega} = \begin{bmatrix} {\bf A}{\bf u}^n + {\bf f}\left({\bf u}^n\right) & {\bf 0} \end{bmatrix} \in {\cal R}^{n\times 2}$. As it is shown in \cite{al2011computing}, the action of $\exp\left(\tau {\bf B}\right)$ on ${\bf v}_1=\begin{bmatrix} 0 & 0 \cdots 0 & 1 & 0\end{bmatrix}^T \in {\cal M}^{(n+2)\times 1}$, 
namely $\exp\left(\tau {\bf B}\right) {\bf v}_1 = \begin{bmatrix} \tau \psi_1\left(\tau {\bf A}\right) \left({\bf A} {\bf u}^{n+1} + {\bf f}\left({\bf u}^n\right)\right) \\ 0 \\ 0 \end{bmatrix}$. 
In other words, to compute our exponential integrator, we need to calculate the action of this "extended" exponential $\exp\left(\tau {\bf B}\right)$ on vector ${\bf v}_1$.
Due to the stability issue, we express
\begin{equation}
\exp\left(\tau {\bf B}\right) {\bf v}_1= \left(\exp \left(\frac{1}{s} \tau {\bf B} \right)\right)^s {\bf v}_1 = 
\exp \left(\frac{1}{s} \tau {\bf B} \right) \cdots \exp \left(\frac{1}{s} \tau {\bf B} \right) {\bf v}_1.
\end{equation}
From Taylor expansion
\begin{equation}
\exp\left(\frac{\tau {\bf B}}{s}\right) = \sum_{j=0}^m \frac{\left(\frac{1}{s}  \tau {\bf B} \right)^j}{j!}.
\end{equation}
We denote
\begin{equation}
{\bf P}_j =\frac{\left(\frac{1}{s}  \tau {\bf B} \right)^j}{j!}.
\end{equation}
With proper selection of $s$ and $m$, we have,
\begin{equation}
\begin{aligned}
\exp\left(\tau {\bf B}\right) {\bf v}_1 & = 
 & \underbrace{\left({\bf P}_0 + {\bf P}_1 + \cdots + {\bf P}_m \right)  \times  \cdots \times \left({\bf P}_0 + {\bf P}_1 + \cdots + {\bf P}_m \right)  }_{s \textrm{ times}}
\end{aligned}.
\end{equation}
From the multinomial theorem, \cite{multinomial}, we have
\begin{equation}
\begin{aligned}
\exp\left(\tau {\bf B}\right) {\bf v}_1 = \sum_{\substack{k_0+k_1+k_2+\cdots+k_m=s\\ k_0,k_1,\cdots,k_m\geq 0}} {{s}\choose{k_1,k_2,\cdots,k_m}} \prod_{t=0}^m {\bf P}_t^{k_t} {\bf v}_1
\end{aligned}.
\end{equation}
This method of computing the action of the exponential integrator allows precomputing the particular powers of ${\bf P}$, from 1 to $s$, and sum them up with proper coefficients.}

\subsection{\revision{Unconditionally stable Crank-Nicolson method}}

\revision{To compare the exponential integrators method with the alternative state-of-the-art method, we describe the Crank-Nicolson time marching scheme.
Namely,
\begin{eqnarray}
\begin{aligned}
{\bf u}^{t+1} =& \frac{\tau}{2} \left( {\bf A} {\bf u}^{t+1} + {\bf A} {\bf u}^{t}\right) + {\bf F}{\bf u}^t, 
\end{aligned}
\end{eqnarray}
which results in a system of linear equations 
\begin{eqnarray}
\begin{aligned}
\left({\bf I} + 
\frac{\tau}{2} {\bf A}\right) {\bf u}^{t+1} &= 
\frac{\tau}{2} {\bf A} {\bf u}^{t} + {\bf F}{\bf u}^t, 
\end{aligned}
\end{eqnarray}
to be solved in every time step.}

\section{Numerical results}
\subsection{MRI scan of the human brain}

\begin{figure}[!htb]
  \centering
\includegraphics[width=0.2\textwidth]{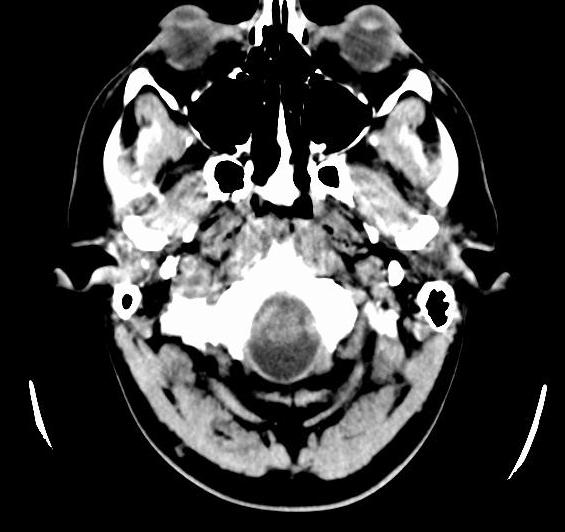}\includegraphics[width=0.2\textwidth] {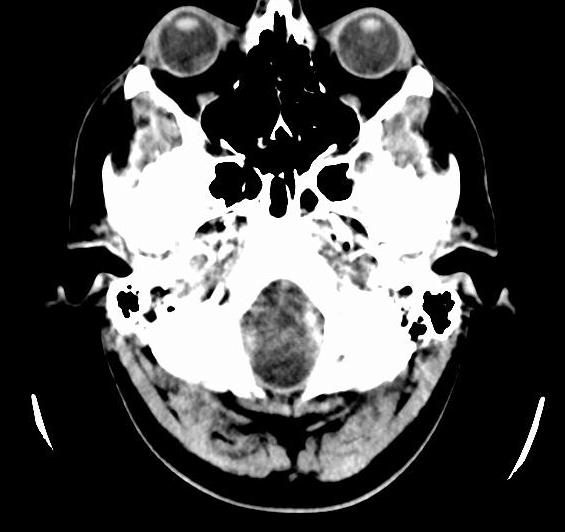}\includegraphics[width=0.2\textwidth]{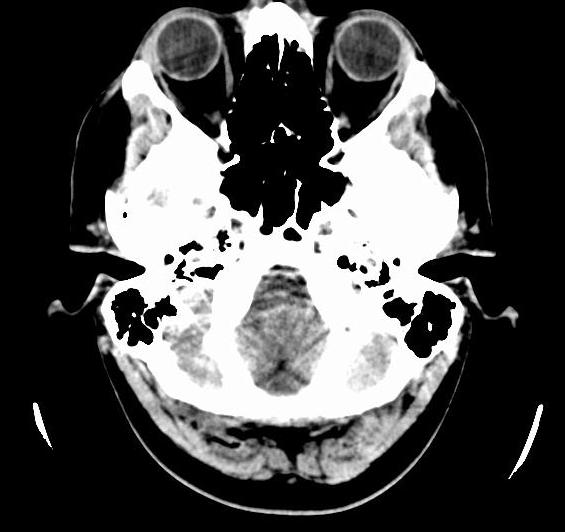}\includegraphics[width=0.2\textwidth]{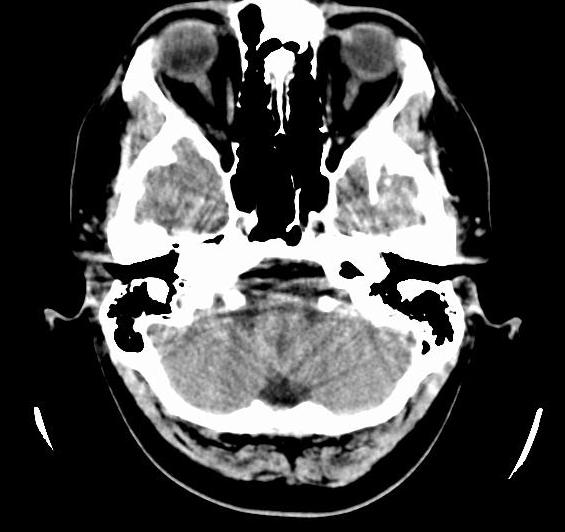}\includegraphics[width=0.2\textwidth]{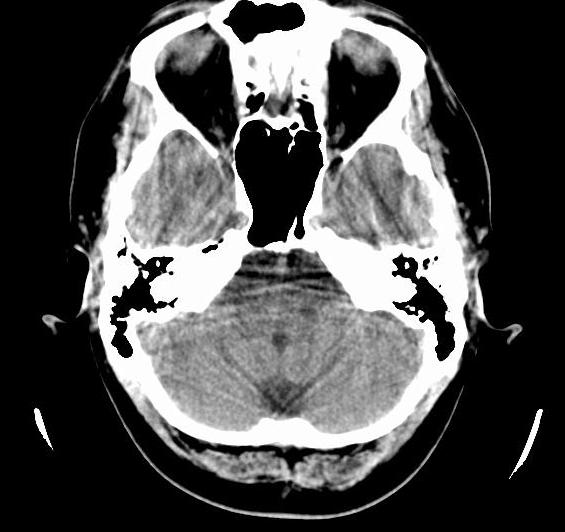}\\ \includegraphics[width=0.2\textwidth]{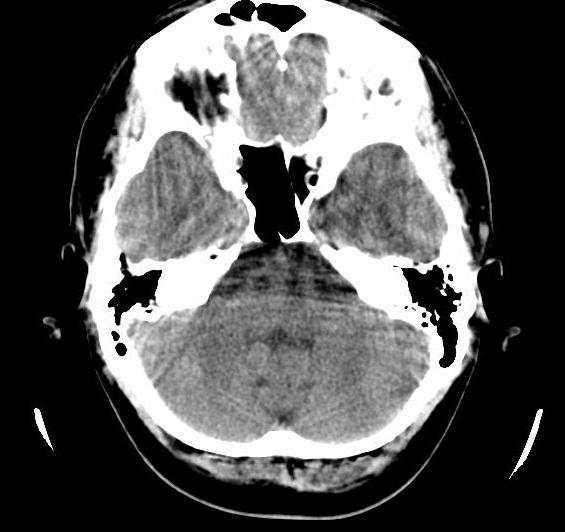}\includegraphics[width=0.2\textwidth] {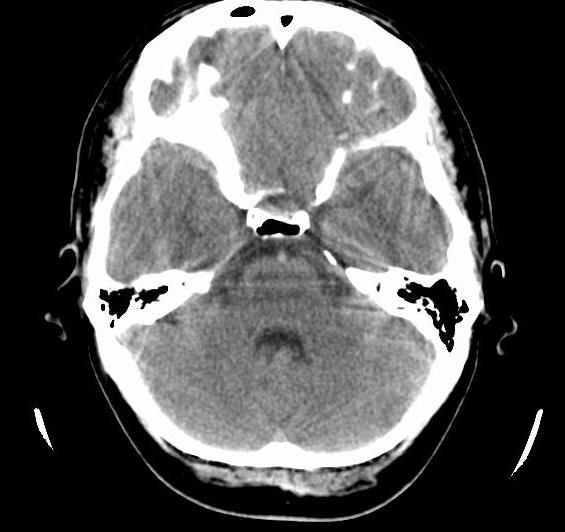}\includegraphics[width=0.2\textwidth]{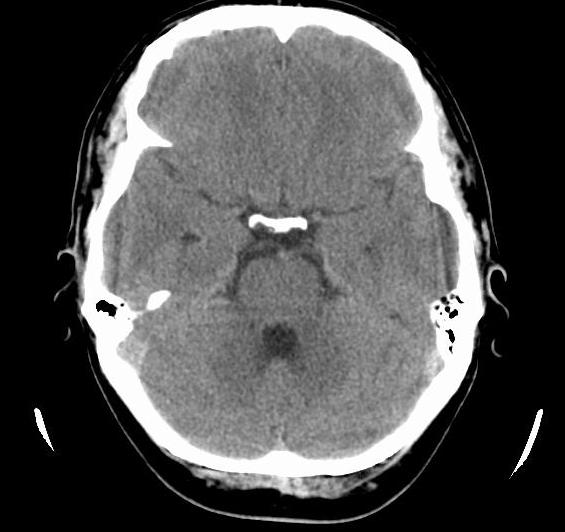}\includegraphics[width=0.2\textwidth]{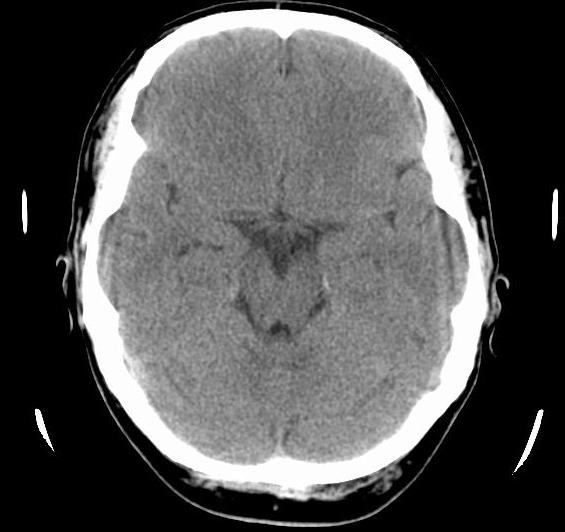}\includegraphics[width=0.2\textwidth]{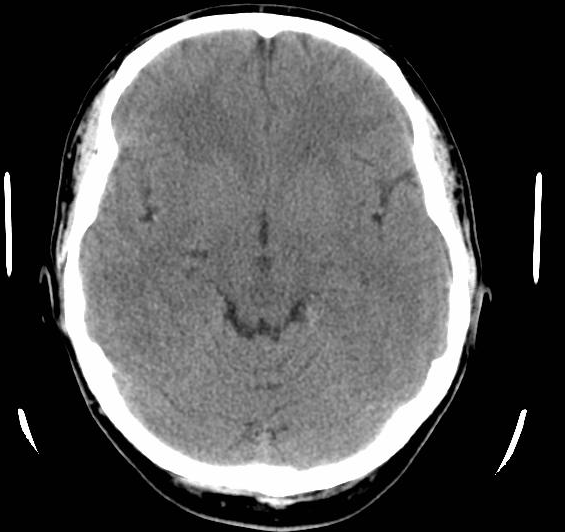}\\ \includegraphics[width=0.2\textwidth]{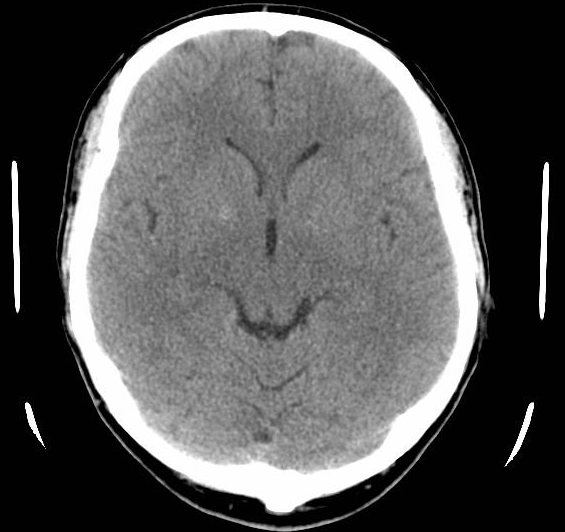}\includegraphics[width=0.2\textwidth] {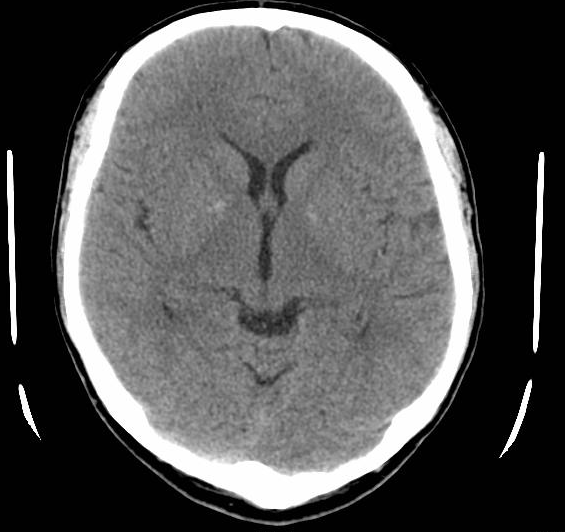}\includegraphics[width=0.2\textwidth]{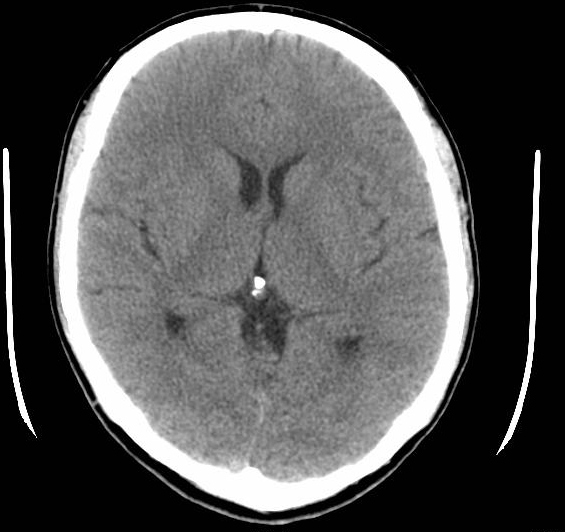}\includegraphics[width=0.2\textwidth]{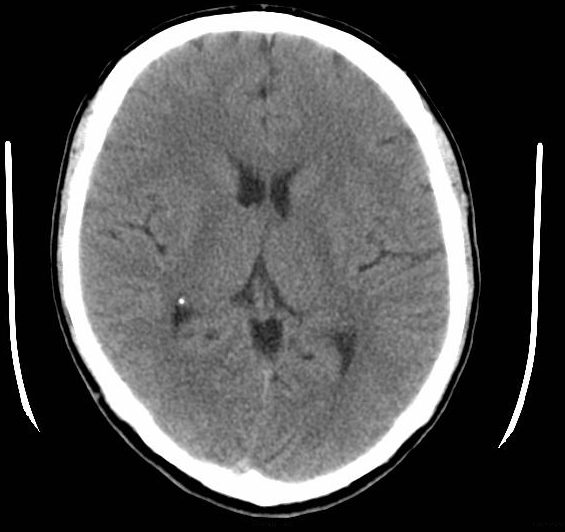}\includegraphics[width=0.2\textwidth]{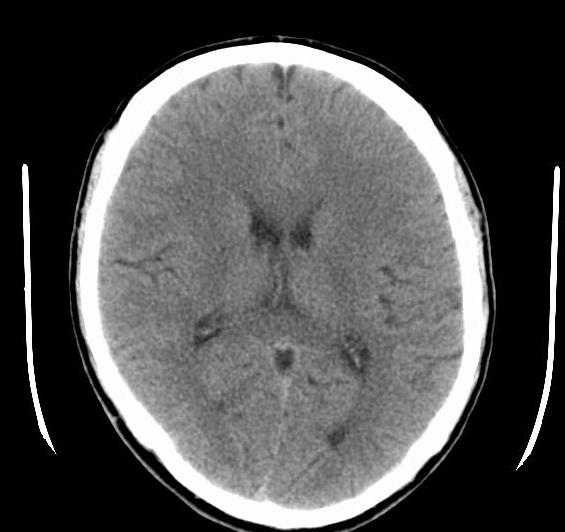}\\ \includegraphics[width=0.2\textwidth]{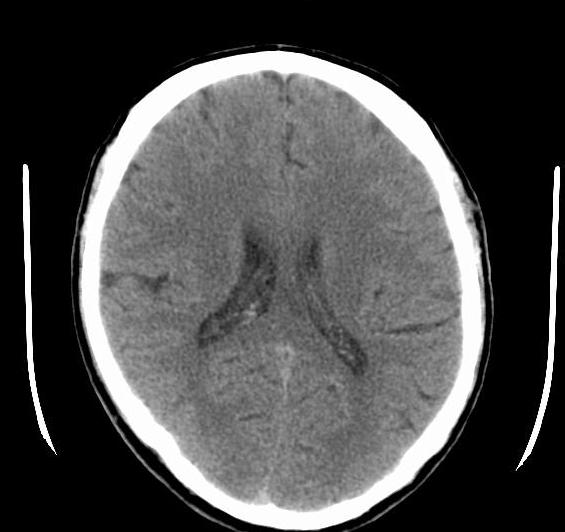}\includegraphics[width=0.2\textwidth] {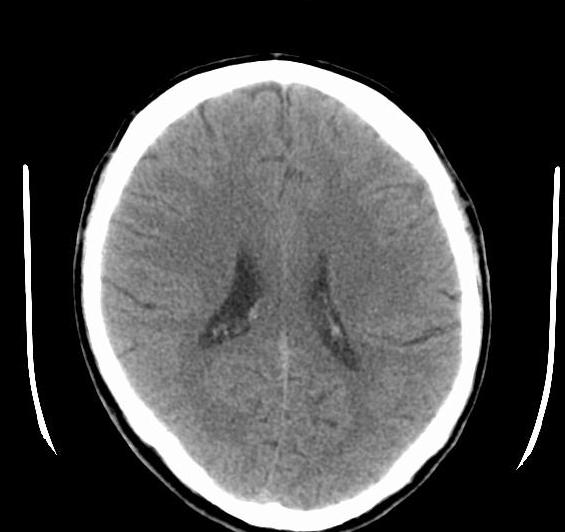}\includegraphics[width=0.2\textwidth]{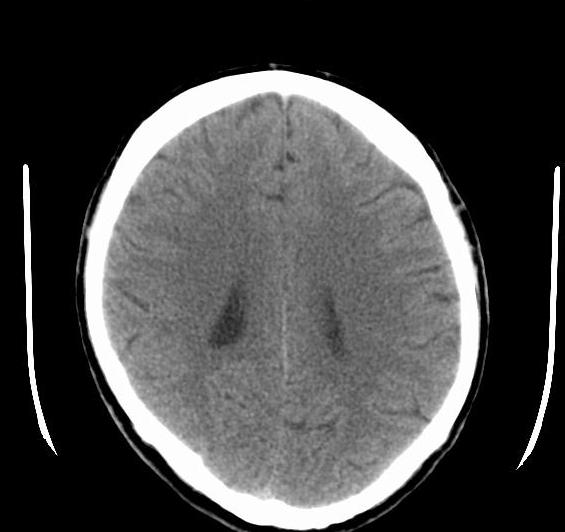}\includegraphics[width=0.2\textwidth]{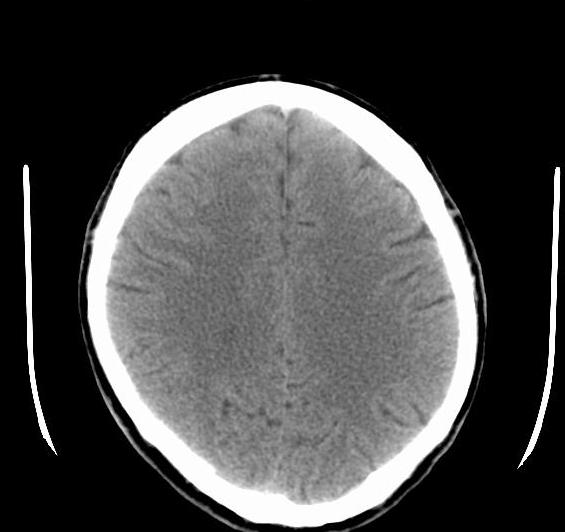}\includegraphics[width=0.2\textwidth]{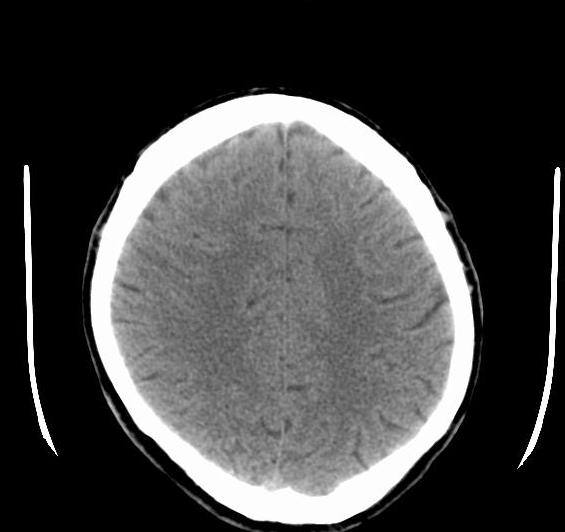}\\  \includegraphics[width=0.2\textwidth]{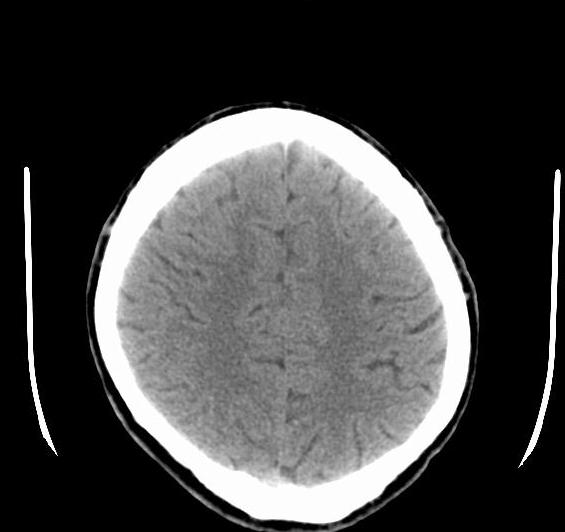}\includegraphics[width=0.2\textwidth] {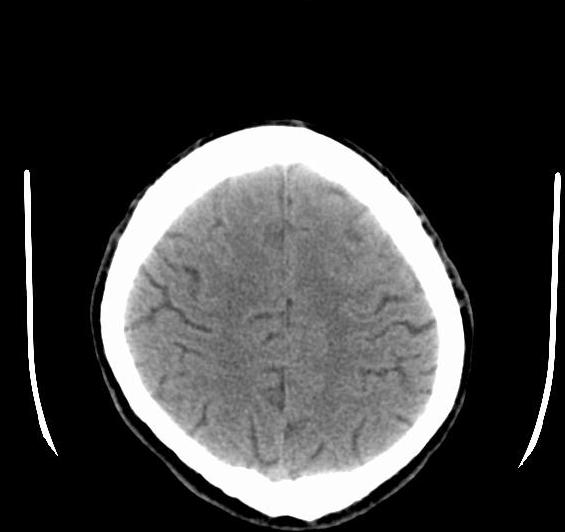}\includegraphics[width=0.2\textwidth]{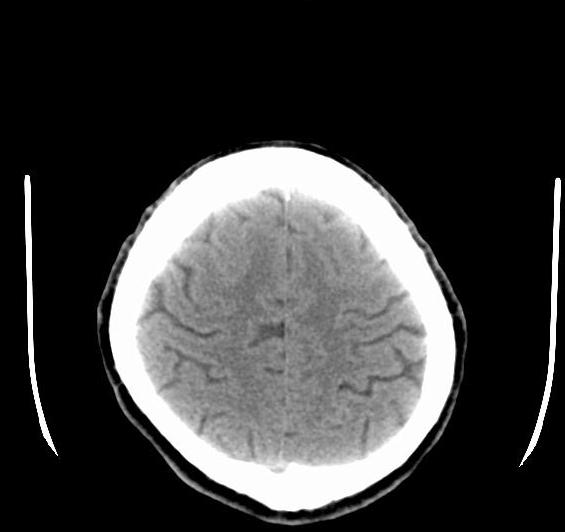}\includegraphics[width=0.2\textwidth]{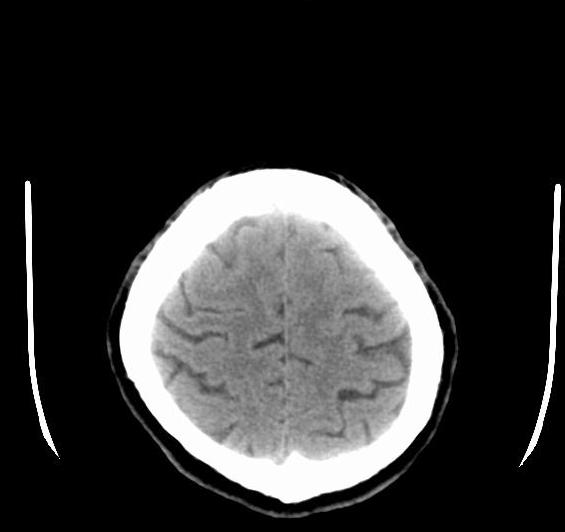}\includegraphics[width=0.2\textwidth]{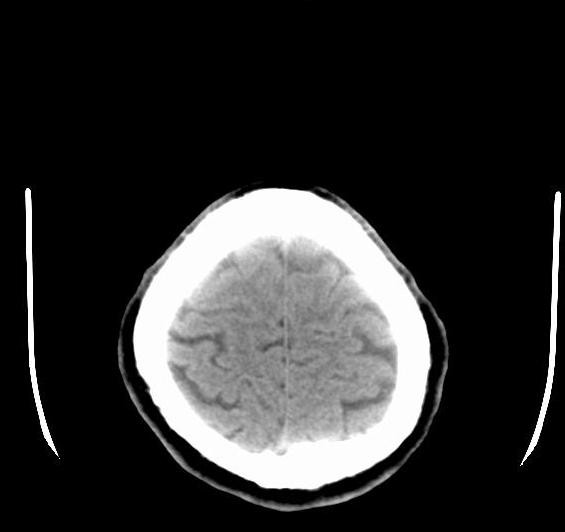}\\  \includegraphics[width=0.2\textwidth]{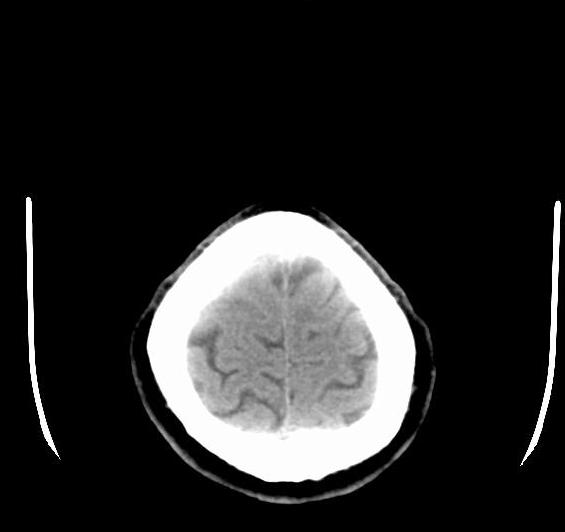}\includegraphics[width=0.2\textwidth] {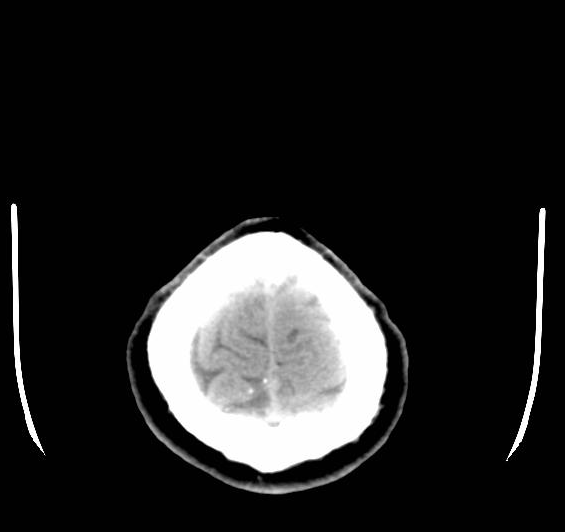}\includegraphics[width=0.2\textwidth]{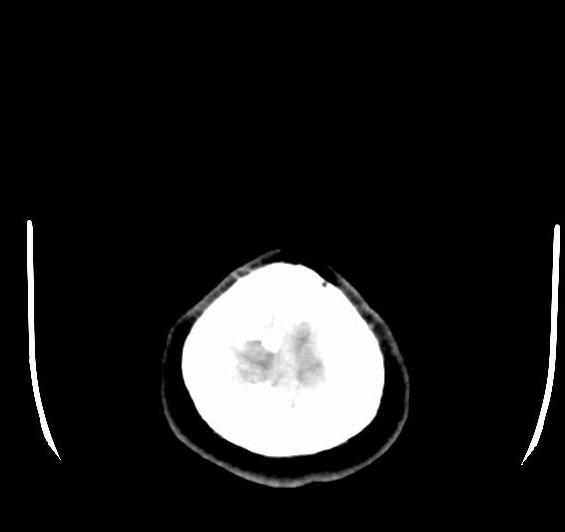}\includegraphics[width=0.2\textwidth]{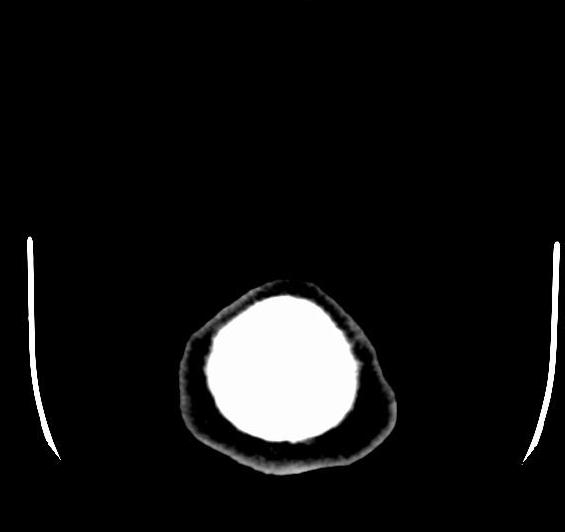}\includegraphics[width=0.2\textwidth]{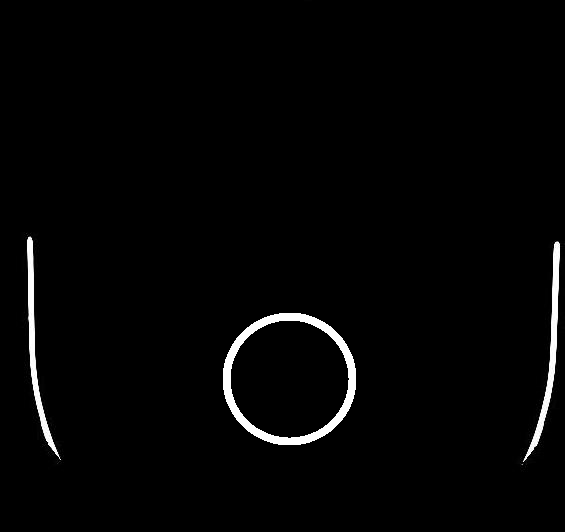}\\   \caption{MRI scans of the head of Maciej Paszy\'nski.}
  \label{fig:head}
\end{figure}
 
Our simulations are based on digital data, and the MRI scan has a sequence of 29 two-dimensional slices, each one with 532 times 565 pixels. Each pixel’s intensity \revision{is a value in [0, 255]}, and it’s proportional to the material’s (skull, skin, white matter or gray matter, and air) normalized density. Exemplary slices of the human head from the MRI scan of the head of one of the authors (Maciej Paszy\'nski) are presented in Figure \ref{fig:head}. Next, according to the MRI scan data, we employ material data changing on the skull, skin, white matter, gray matter, and air. {We assume air (MRI scan data $\leq$ 1), white matter  (1 $\leq$ approximation $\leq$ 115), gray matter (115 $<$ approximation $\leq$ 240) and skull (approximation $>$ 240). }

\subsection{Two-dimensional glioblastoma brain tumor simulation}

We consider one central slice of the MRI scan data.
Let $\Omega\subset [0,200]^2$ [mm$^2$] be an open set and $I=(150,3500)$ [days] the simulated time interval. We solve the two-dimensional version of the Fisher-Kolmogorov diffusion-reaction equation with logistic growth describing the brain tumor dynamics (\ref{PDE}) \cite{Logistic}.
We prescribe the following initial state
\begin{equation}\label{PDE2DBC}
u(x_1,x_2;t_0)=0.1 \exp(-10((x_1-x_1^{init})^2+(x_2-x_2^{init})^2)).
\end{equation}
where $(x_1^{init},x_2^{init})$ is the initial tumor location.
We select $\rho=2.5\times 0.01$ [days$^{-1}$] \cite{c25}. The diffusion coefficient for the white matter is defined as $D_w=1.3\times 0.1$ [mm$^2$days$^{-1}$] \cite{c25}.
The diffusion coefficient for the gray matter is defined as $D_w=0.13\times 0.1$ [mm$^2$days$^{-1}$] \cite{Menze}.

With this parameters setup, we start from $t_0=150$ [days] with $T=3500$, and we perform 100-time steps of the exponential integrator simulation.
Thus, one time step represents $\frac{3500-150}{100}=33$ [days].
The spatial grid for the finite difference method consists of $193 \times 193$ points.
The single bitmap resolution is 532 times 565 pixels.
The 100 time steps exponential integrators method executed on a laptop with Win10, using MATLAB, with 11th Gen Intel(R) Core(TM) i5-11500H @ 2.90GHz, 2.92 GHz, and 32 GB of RAM, takes less than 6 seconds.

\begin{figure}[!htb]
  \centering
         \includegraphics[width=\textwidth]{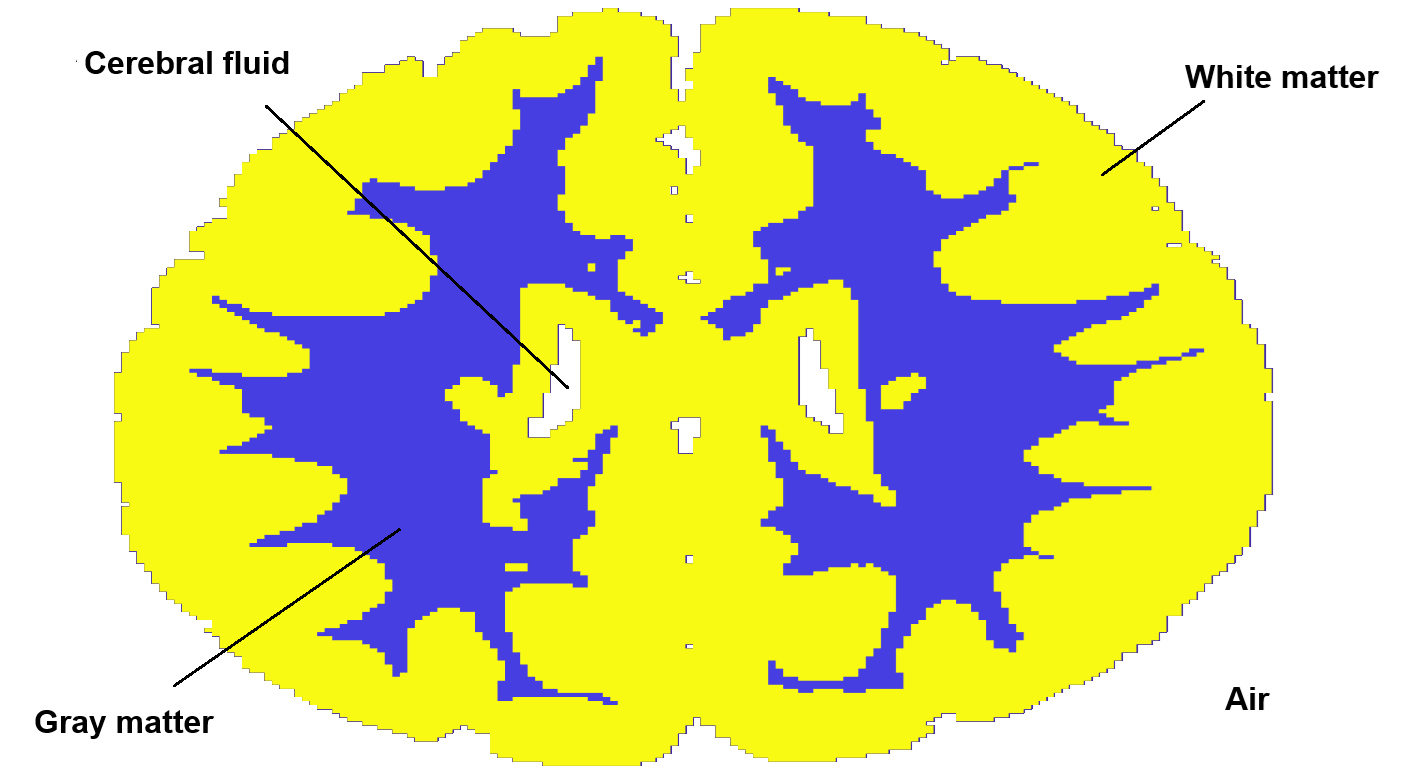}
         \caption{\revision{The map of the white and gray matter, cerebral fluid, and the air, on the axial cross-section.}}
         \label{fig:matters}
\end{figure}

\begin{figure}[!htb]
  \centering
     \begin{subfigure}[b]{0.49\textwidth}
         \centering
         \includegraphics[width=\textwidth]{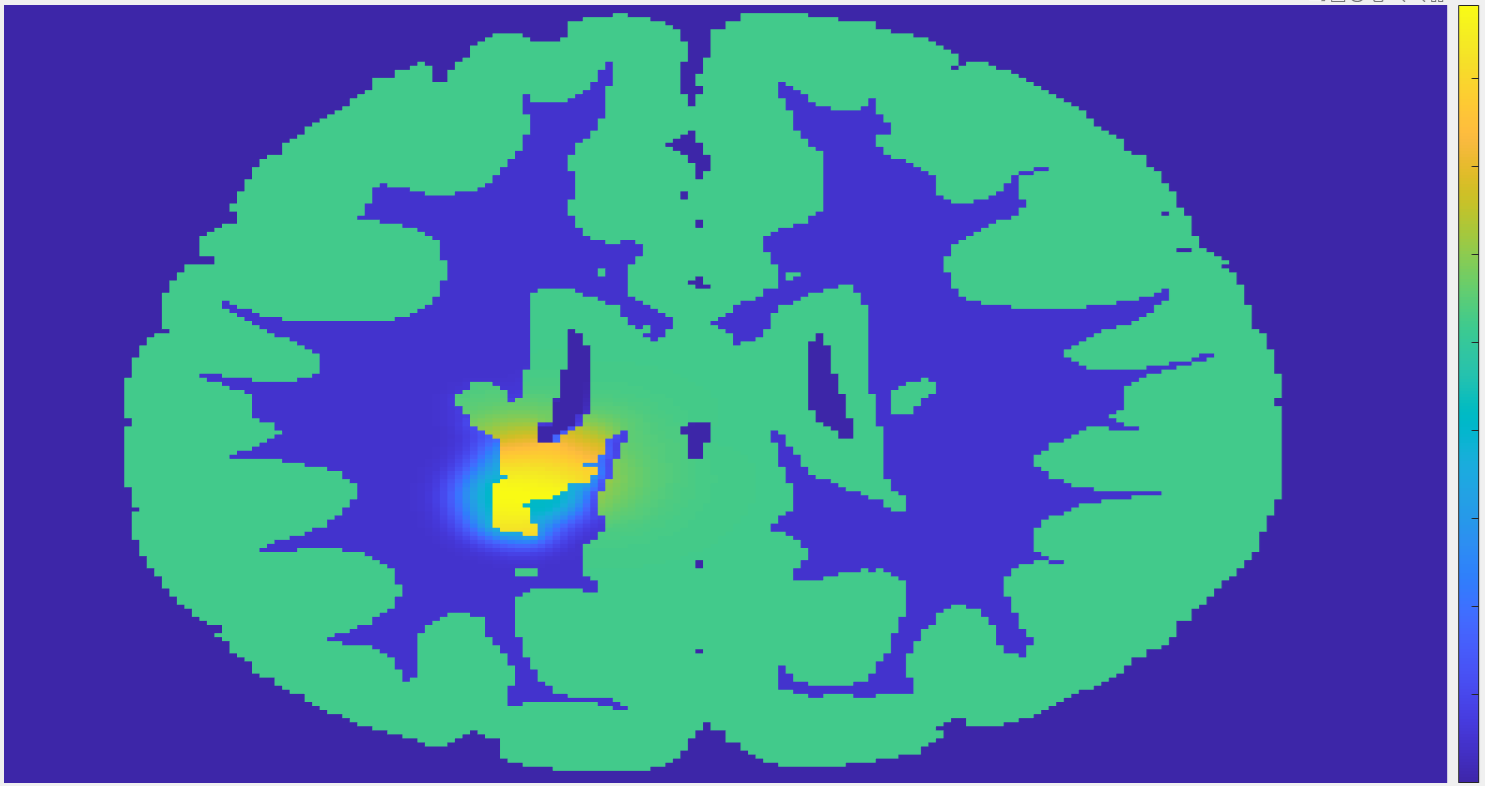}
         \caption{Snapshot at day 418.}
         \label{fig:k}
     \end{subfigure}
     \hfill
     \begin{subfigure}[b]{0.49\textwidth}
         \centering         \includegraphics[width=\textwidth]{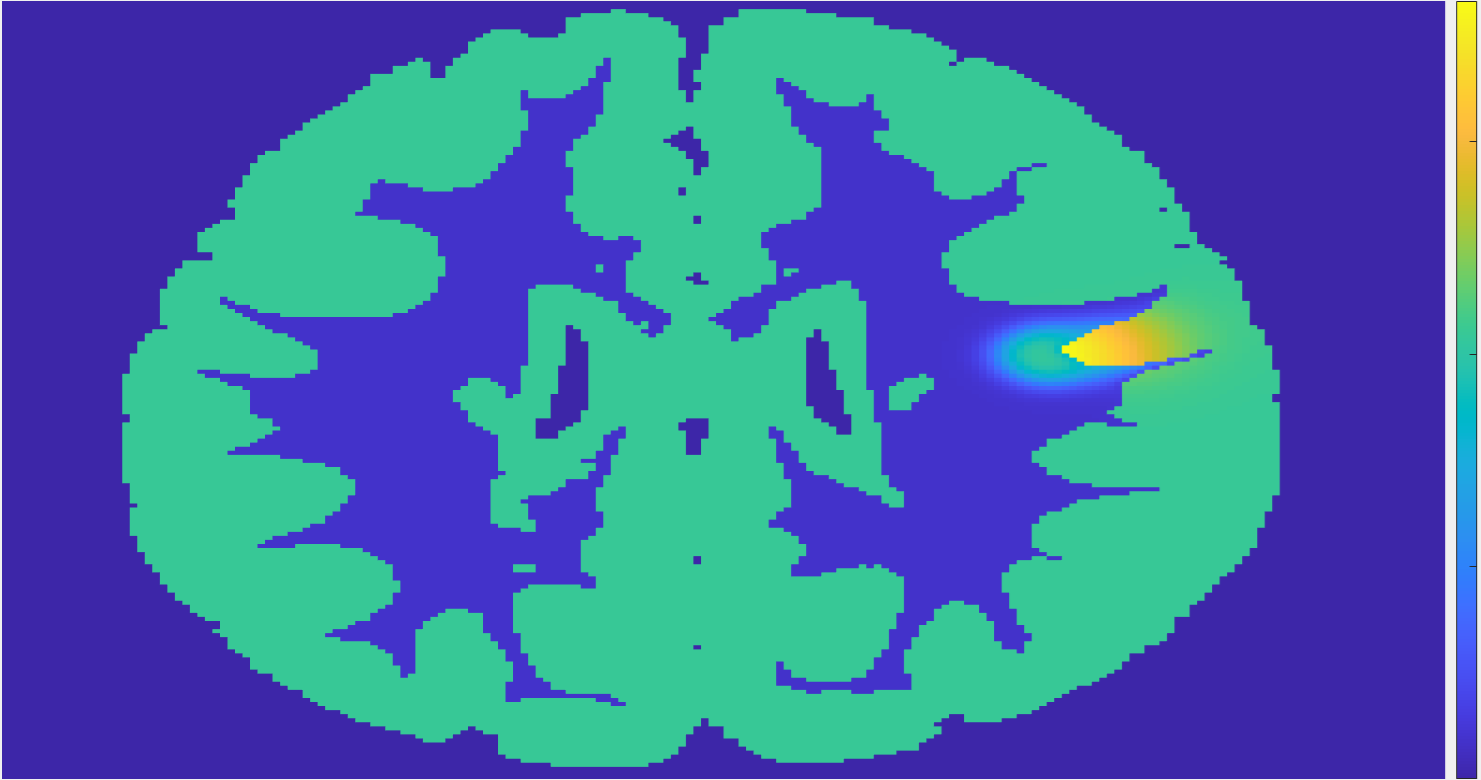}
         \caption{Snapshot at day 418.}
         \label{fig:j}
     \end{subfigure}
\\     
     \begin{subfigure}[b]{0.49\textwidth}
         \centering
         \includegraphics[width=\textwidth]{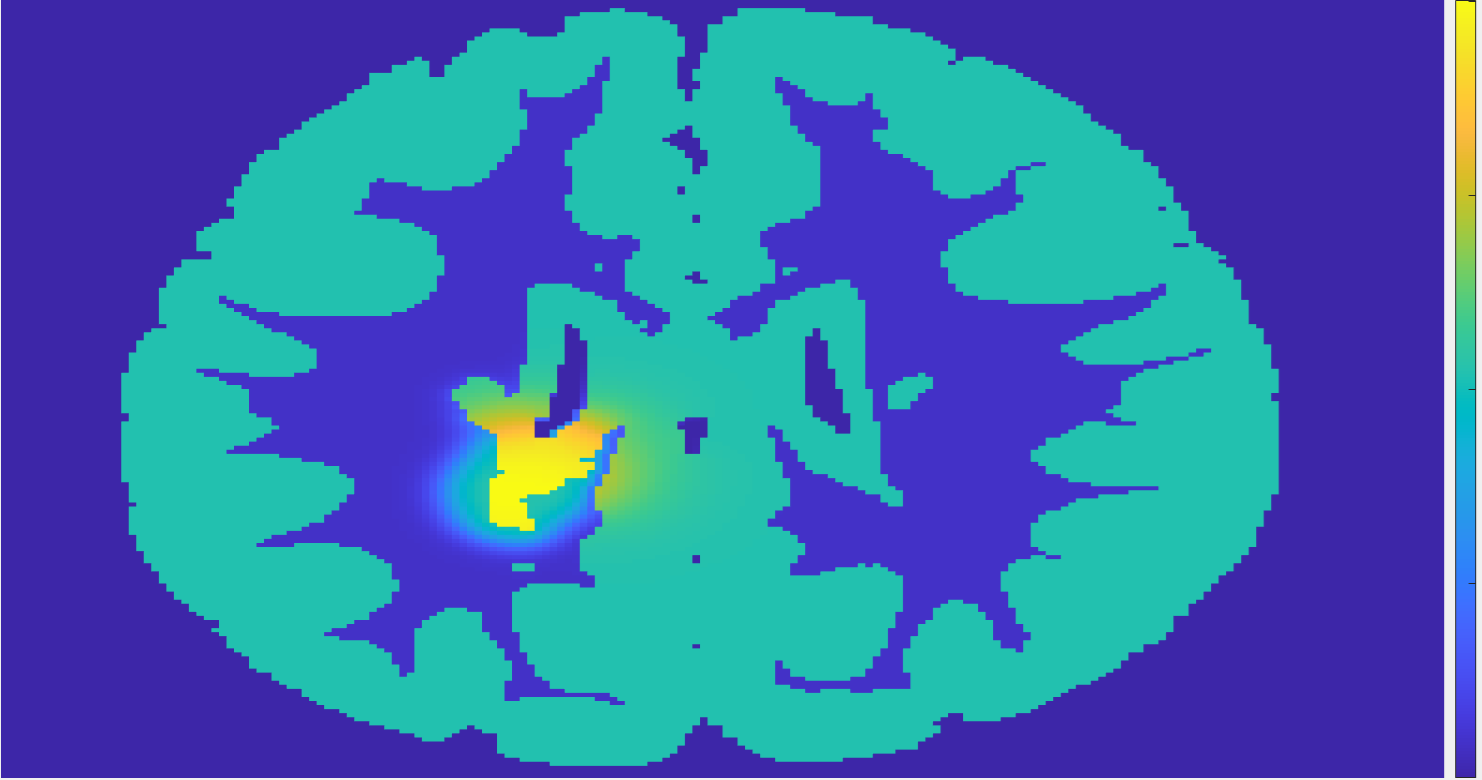}
         \caption{Snapshot at day 552.} 
         \label{fig:m}
     \end{subfigure}
     \hfill
     \begin{subfigure}[b]{0.49\textwidth}
         \centering         \includegraphics[width=\textwidth]{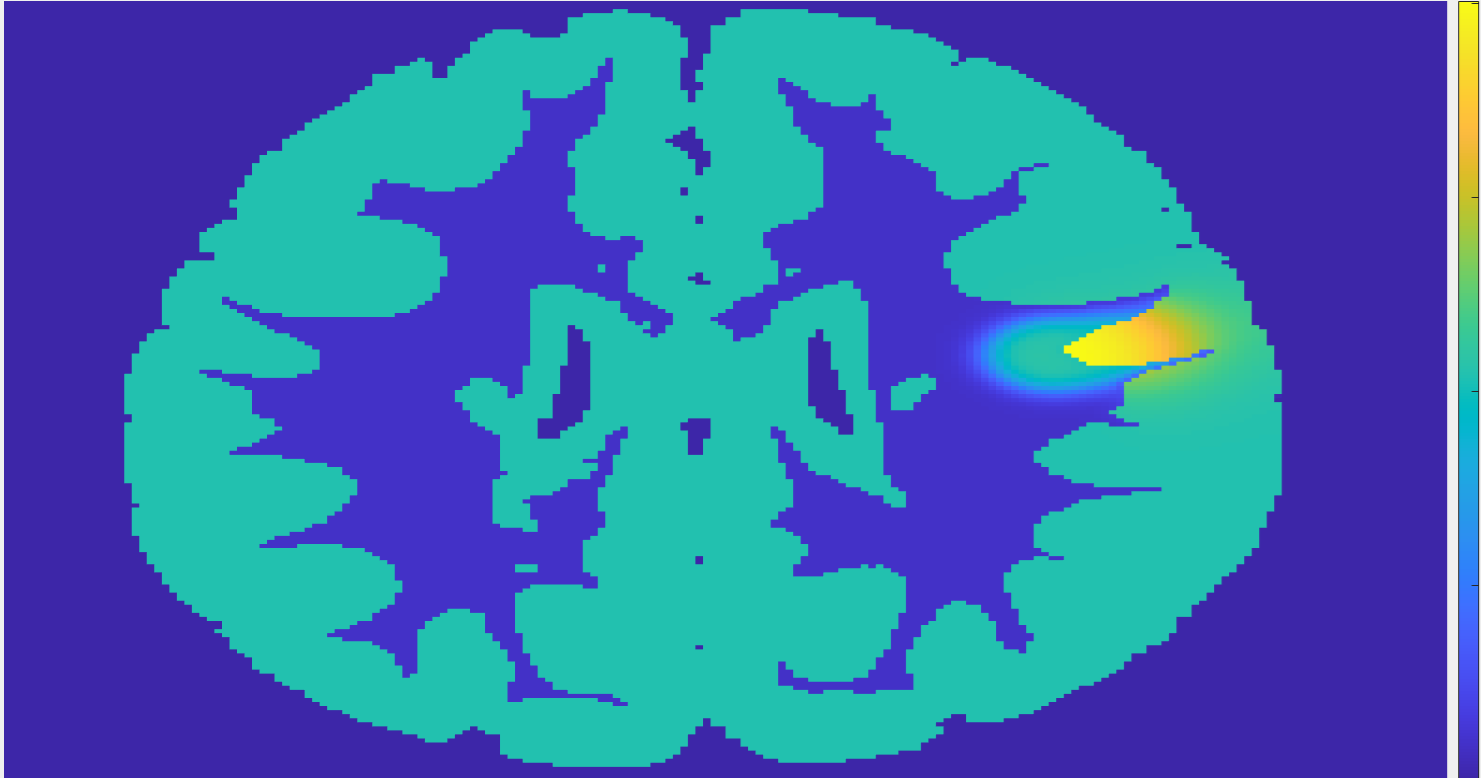}
         \caption{Snapshot at day 552.} 
         \label{fig:n}
     \end{subfigure}
\\ 
     \begin{subfigure}[b]{0.49\textwidth}
         \centering
         \includegraphics[width=\textwidth]{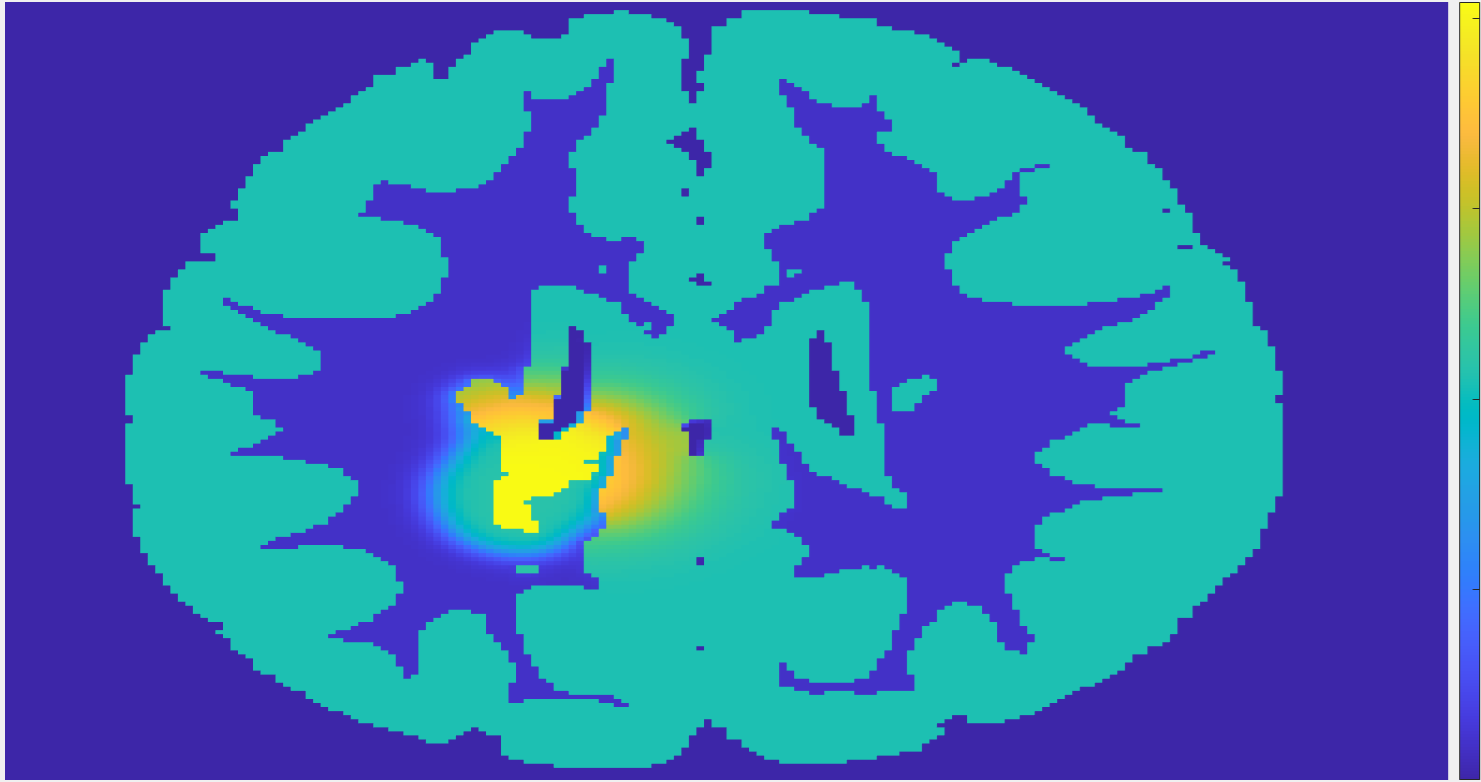}
         \caption{Snapshot at day 686}
         \label{fig:q}
     \end{subfigure}
     \hfill
     \begin{subfigure}[b]{0.49\textwidth}
         \centering         \includegraphics[width=\textwidth]{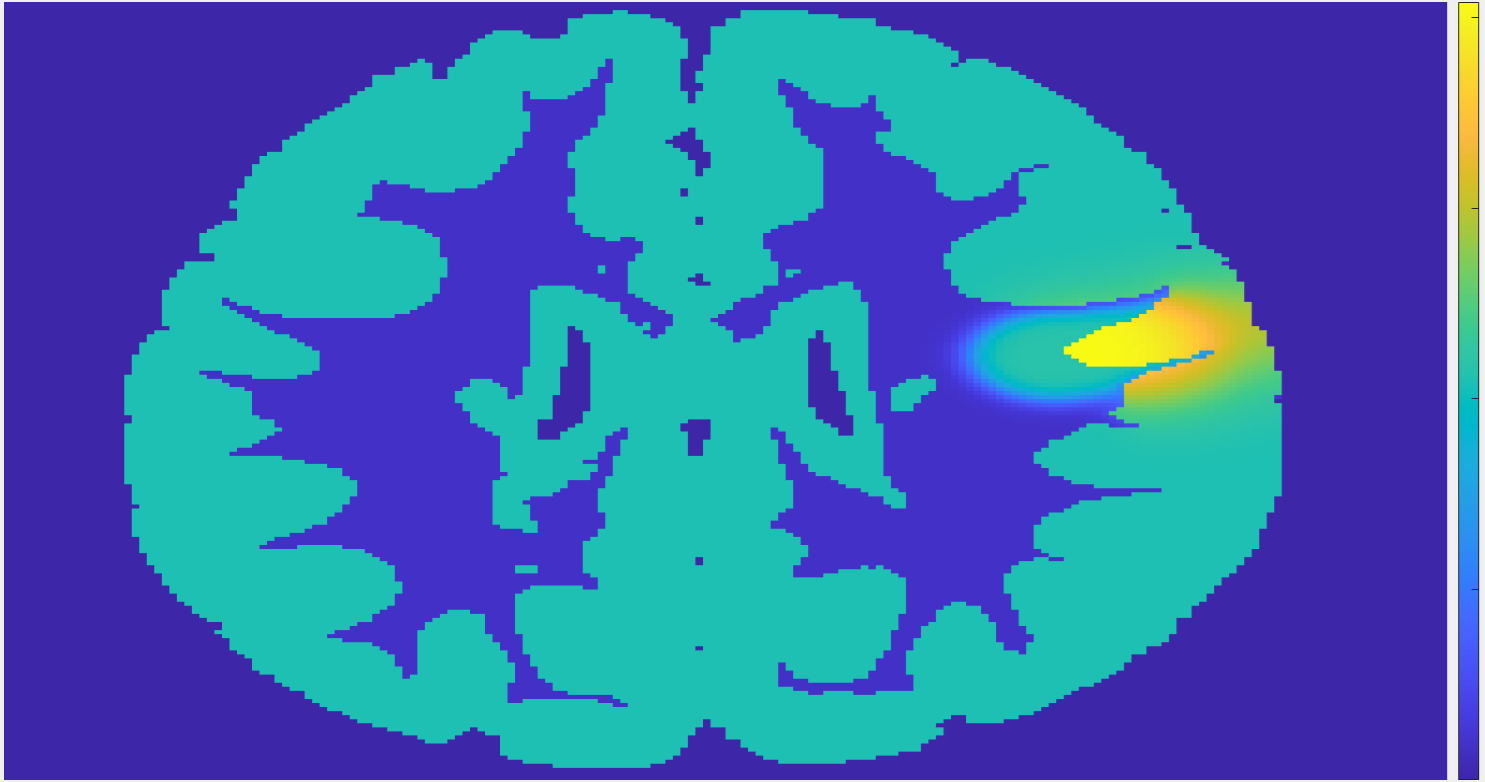}
         \caption{Snapshot at day 686.}
         \label{fig:r} 
     \end{subfigure}
\\
     \begin{subfigure}[b]{0.49\textwidth}
         \centering
         \includegraphics[width=\textwidth]{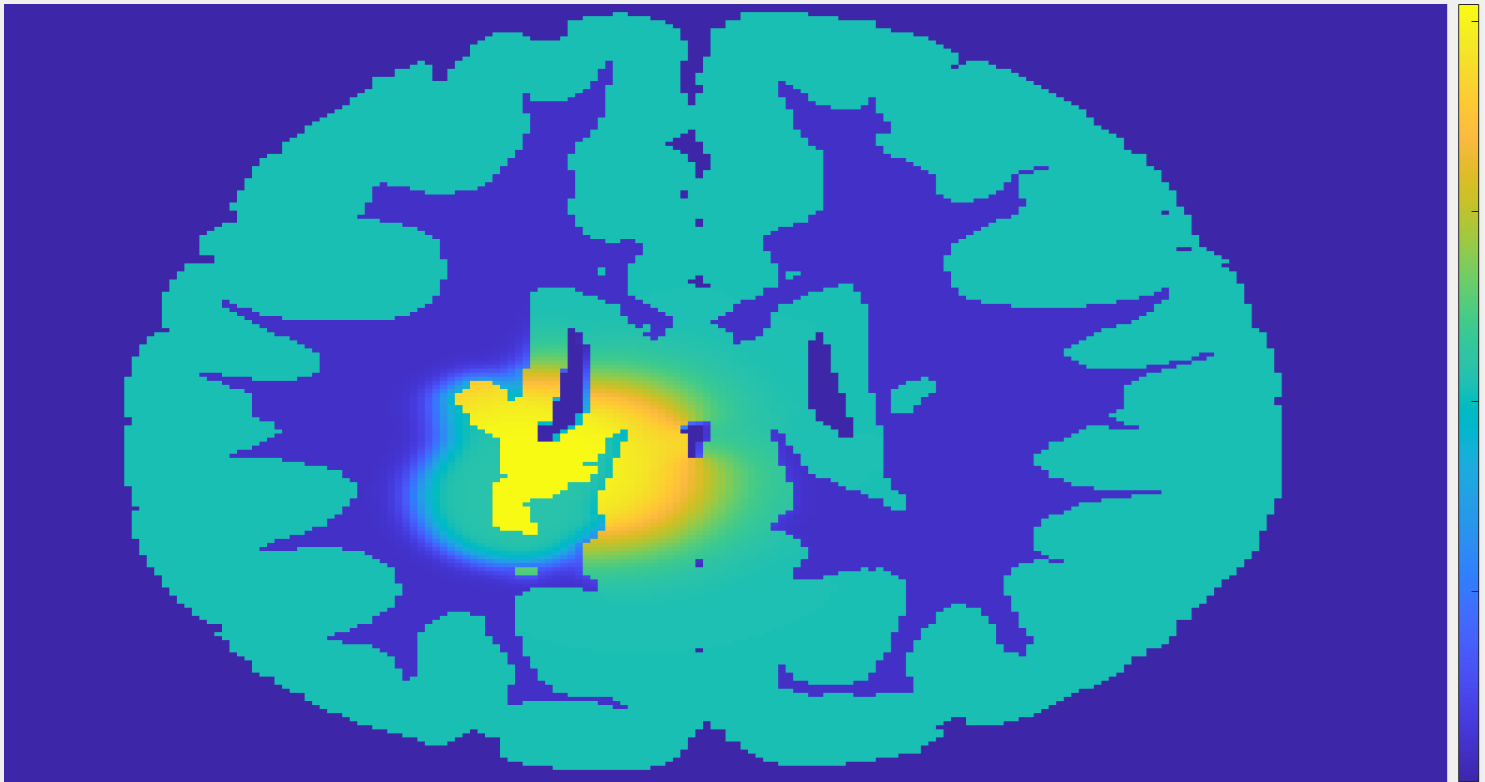}
         \caption{Snapshot at day 820.}
         \label{fig:q}
     \end{subfigure}
     \hfill
     \begin{subfigure}[b]{0.49\textwidth}
         \centering         \includegraphics[width=\textwidth]{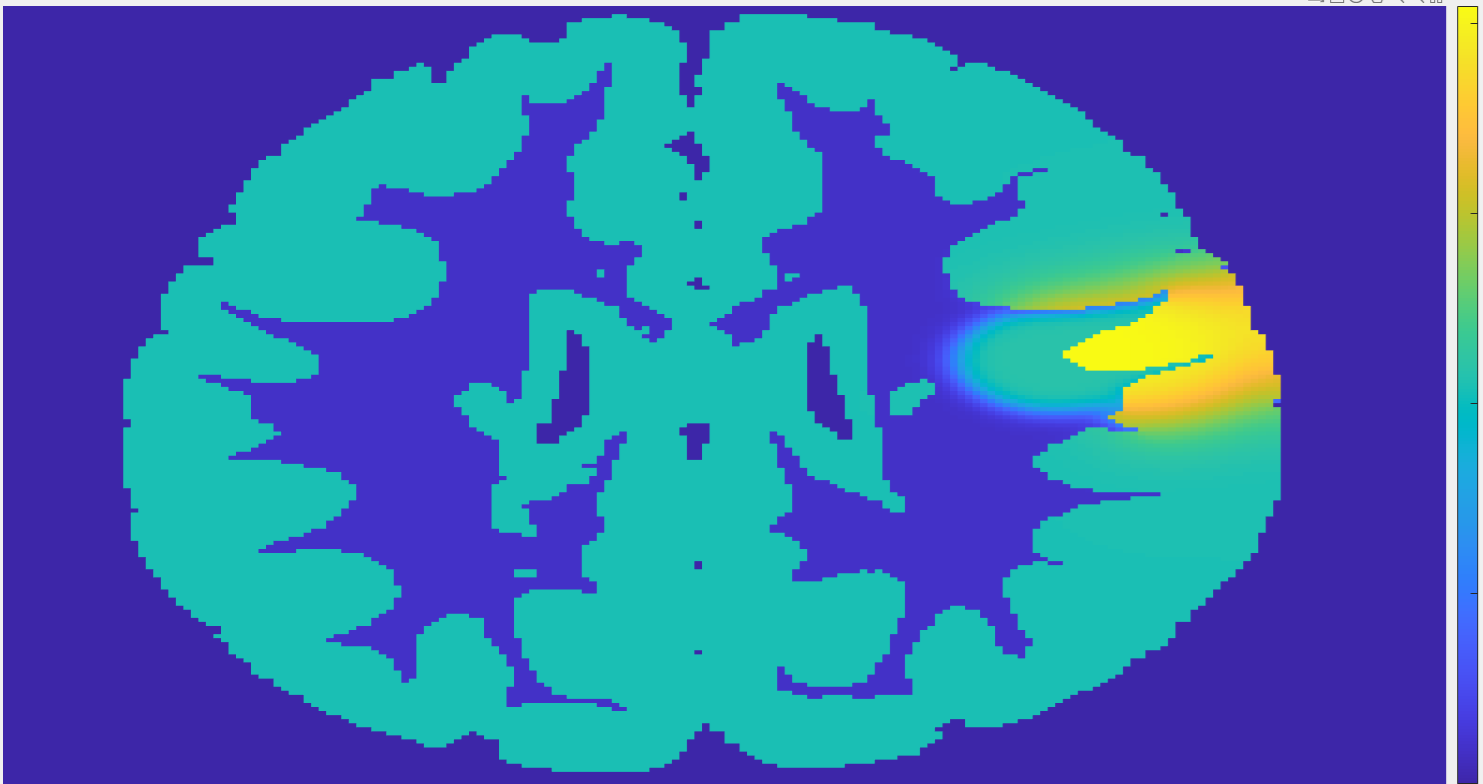}
         \caption{Snapshot at day 820.}
         \label{fig:r} 
     \end{subfigure}
\caption{\revision{Comparison of tumor growth starting in two different brain locations. Snapshots illustrate different growth rates in white and gray matter.}}
  \label{fig:comparison1}
\end{figure}

\begin{figure}[!htb]
  \centering
     \begin{subfigure}[b]{0.49\textwidth}
         \centering
         \includegraphics[width=\textwidth]{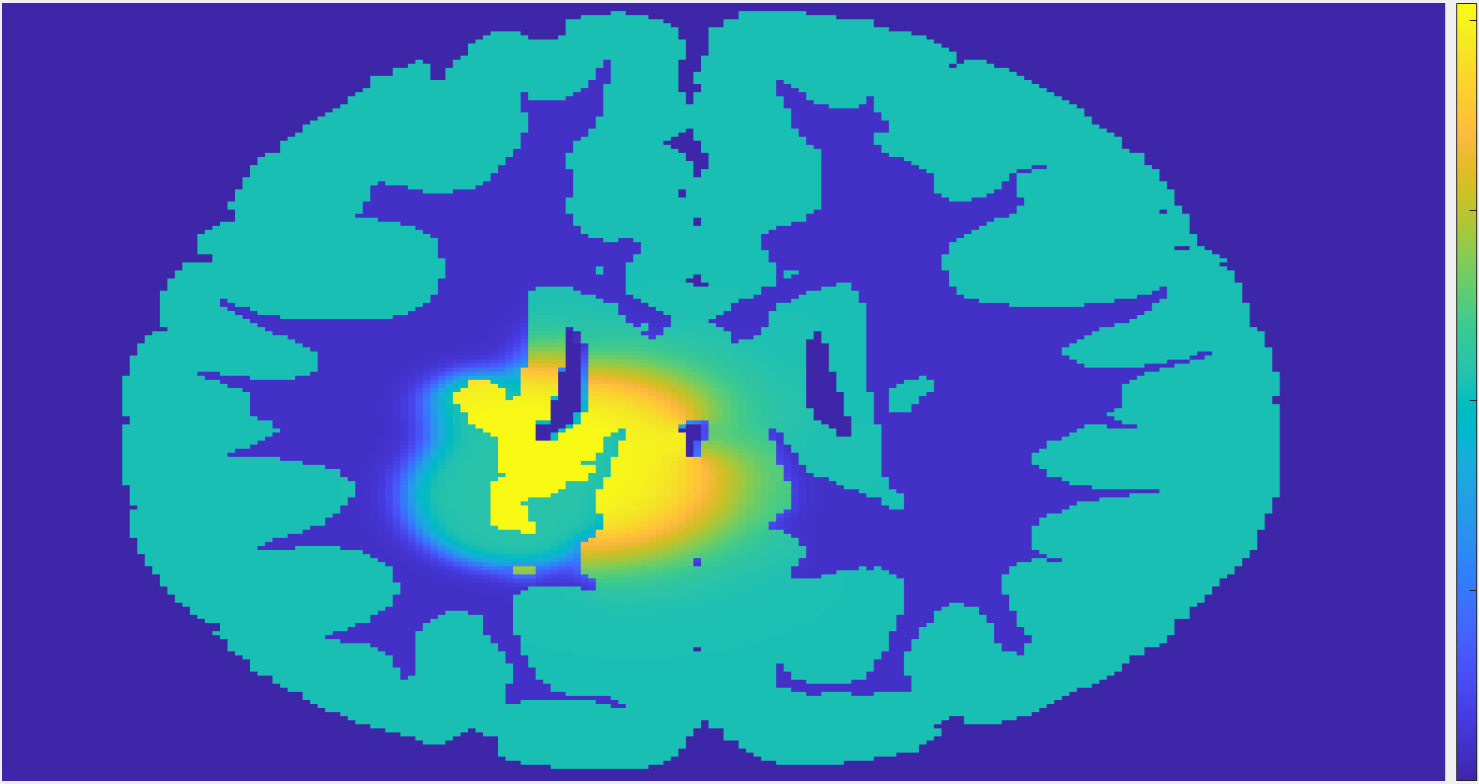}
         \caption{Snapshot at day 954}
         \label{fig:q}
     \end{subfigure}
     \hfill
     \begin{subfigure}[b]{0.49\textwidth}
         \centering         \includegraphics[width=\textwidth]{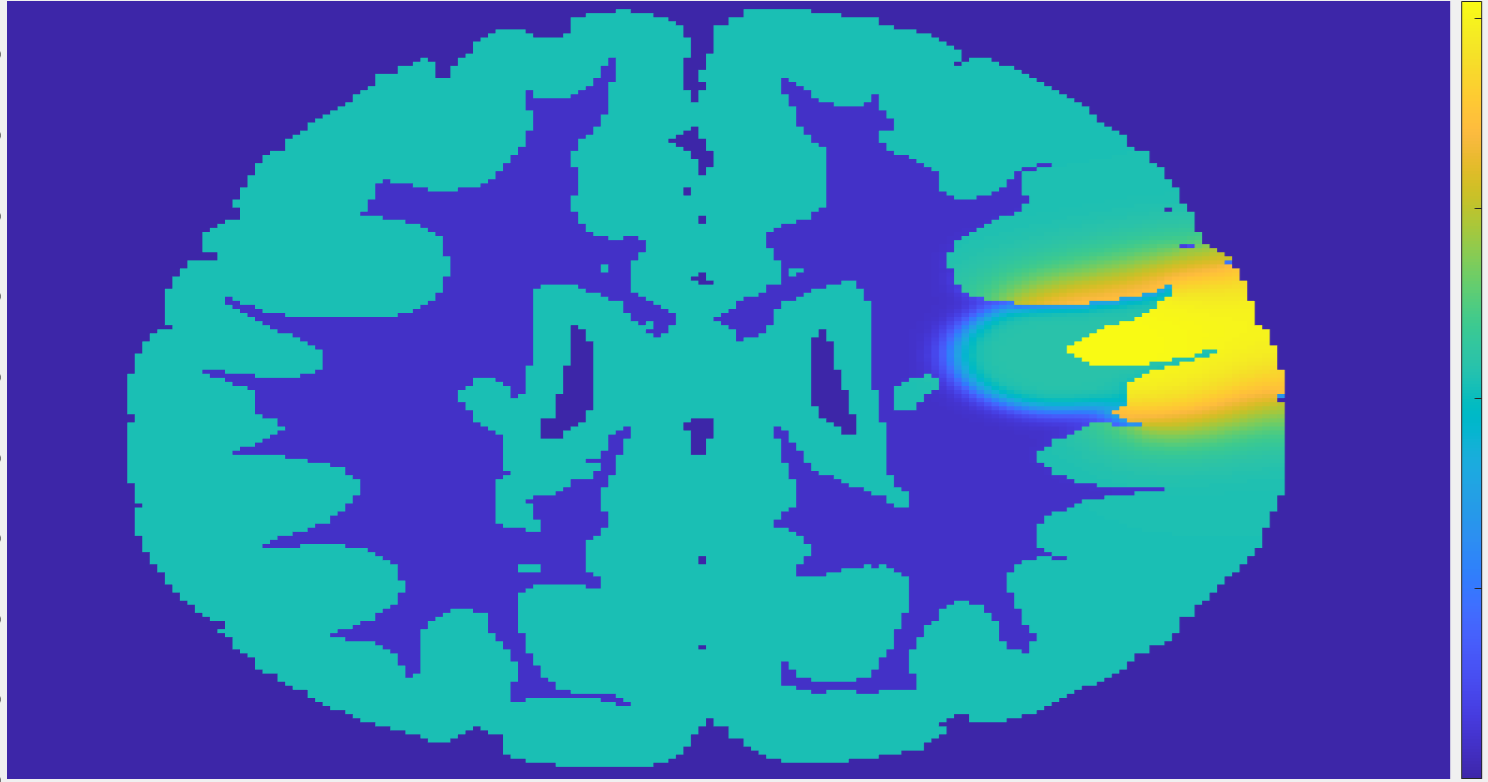}
         \caption{Snapshot at day 954.}
         \label{fig:r} 
     \end{subfigure}
\\
     \begin{subfigure}[b]{0.49\textwidth}
         \centering
         \includegraphics[width=\textwidth]{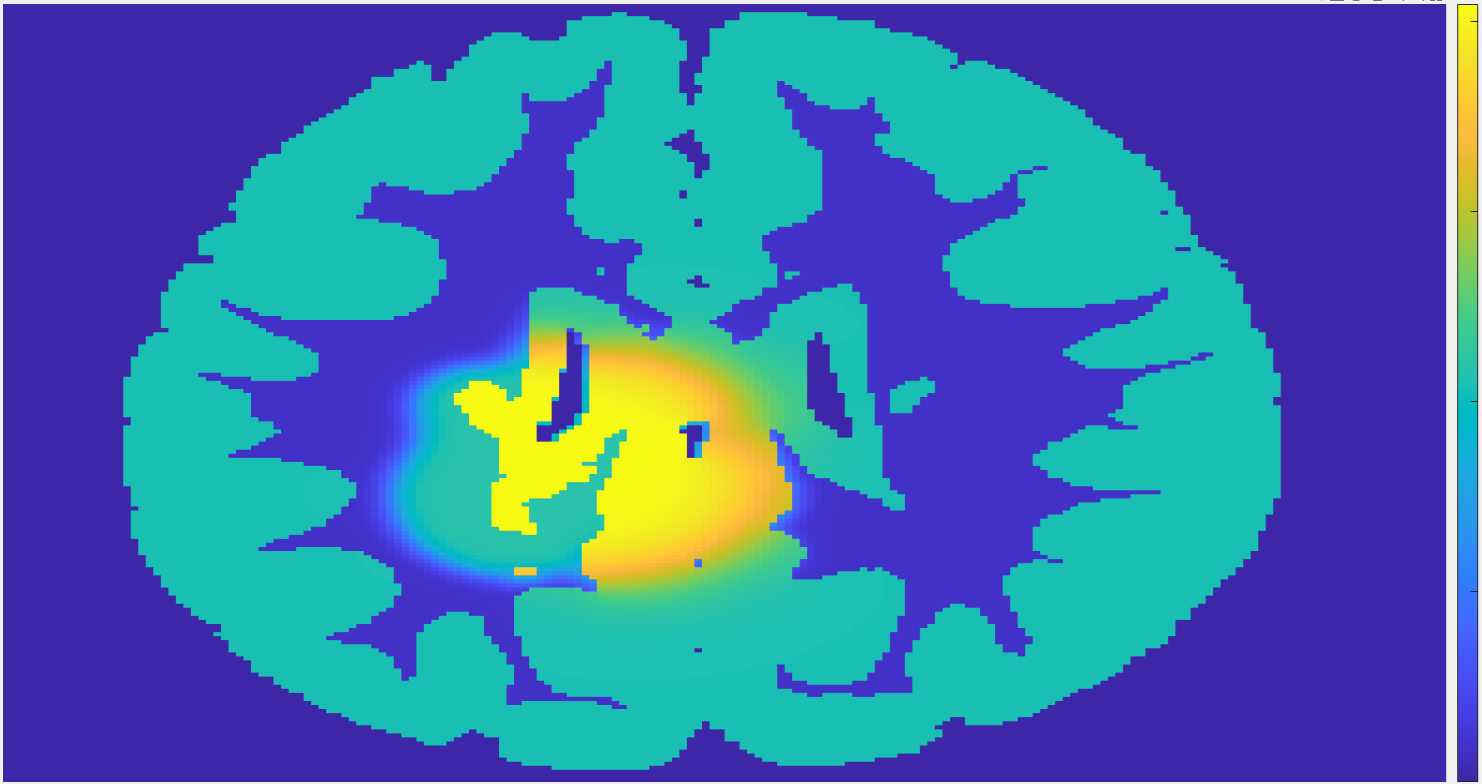}
         \caption{Snapshot at day 1088.}
         \label{fig:q}
     \end{subfigure}
     \hfill
     \begin{subfigure}[b]{0.49\textwidth}
         \centering         \includegraphics[width=\textwidth]{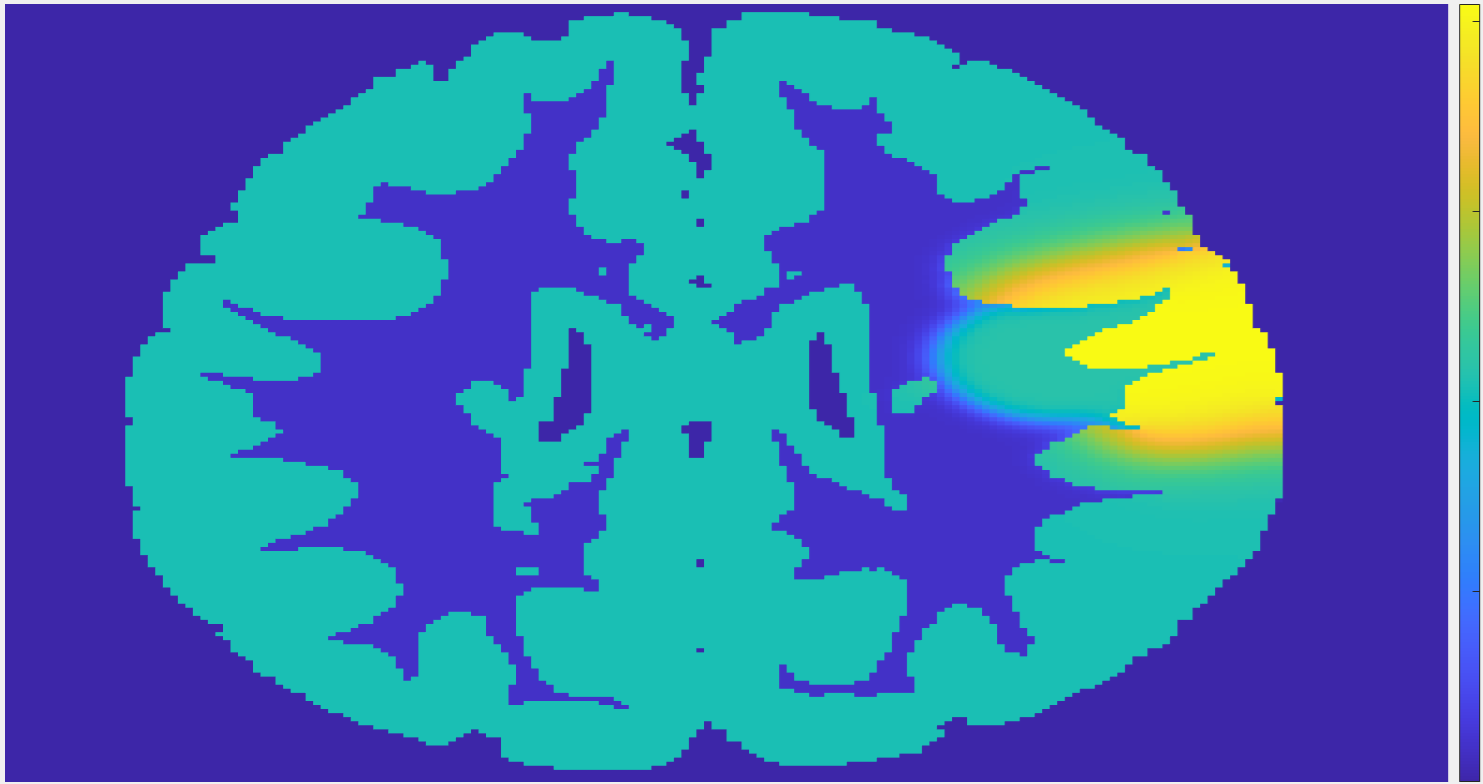}
         \caption{Snapshot at day 1088.}
         \label{fig:r} 
     \end{subfigure}
\\
     \begin{subfigure}[b]{0.49\textwidth}
         \centering
         \includegraphics[width=\textwidth]{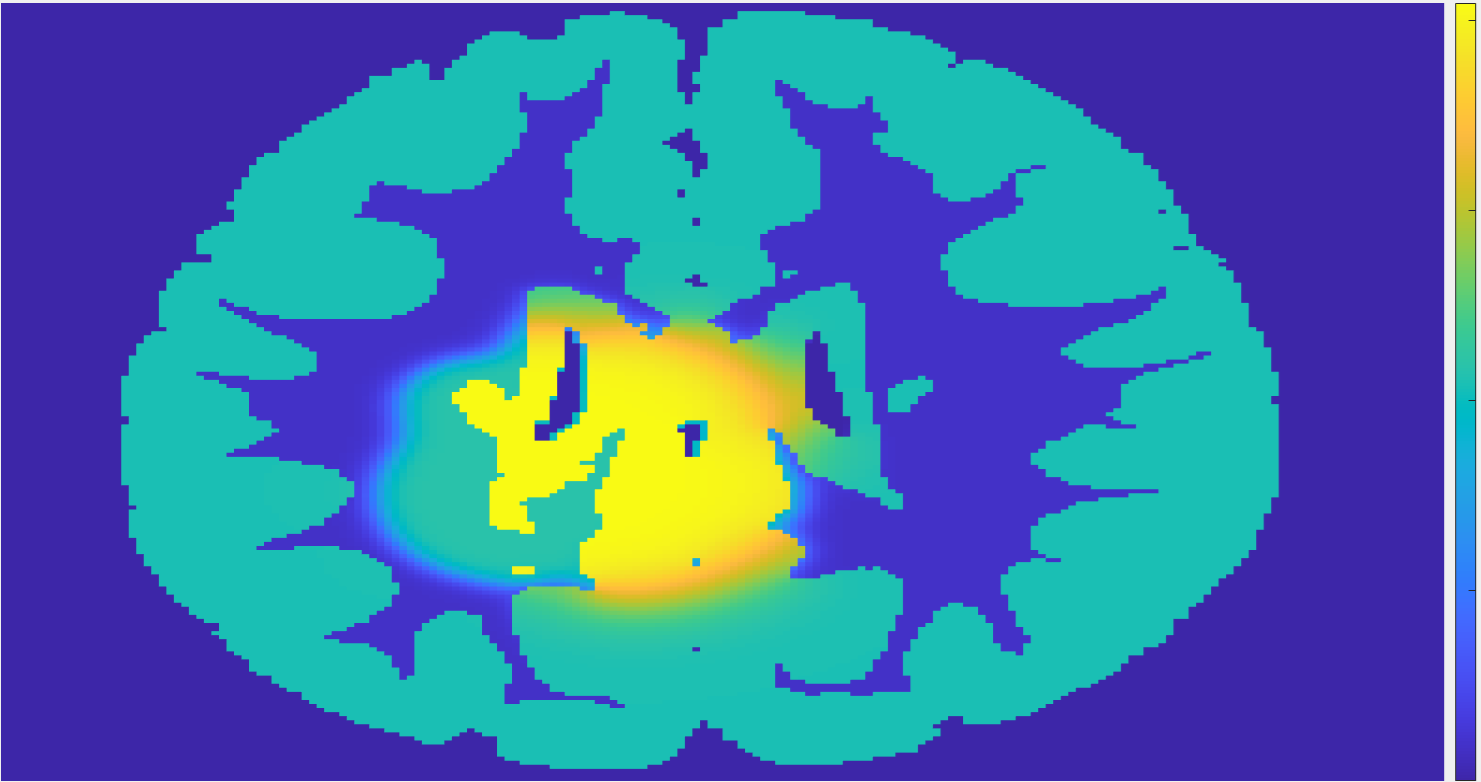}
         \caption{Snapshot at day 1222.}
         \label{fig:q}
     \end{subfigure}
     \hfill
     \begin{subfigure}[b]{0.49\textwidth}
         \centering         \includegraphics[width=\textwidth]{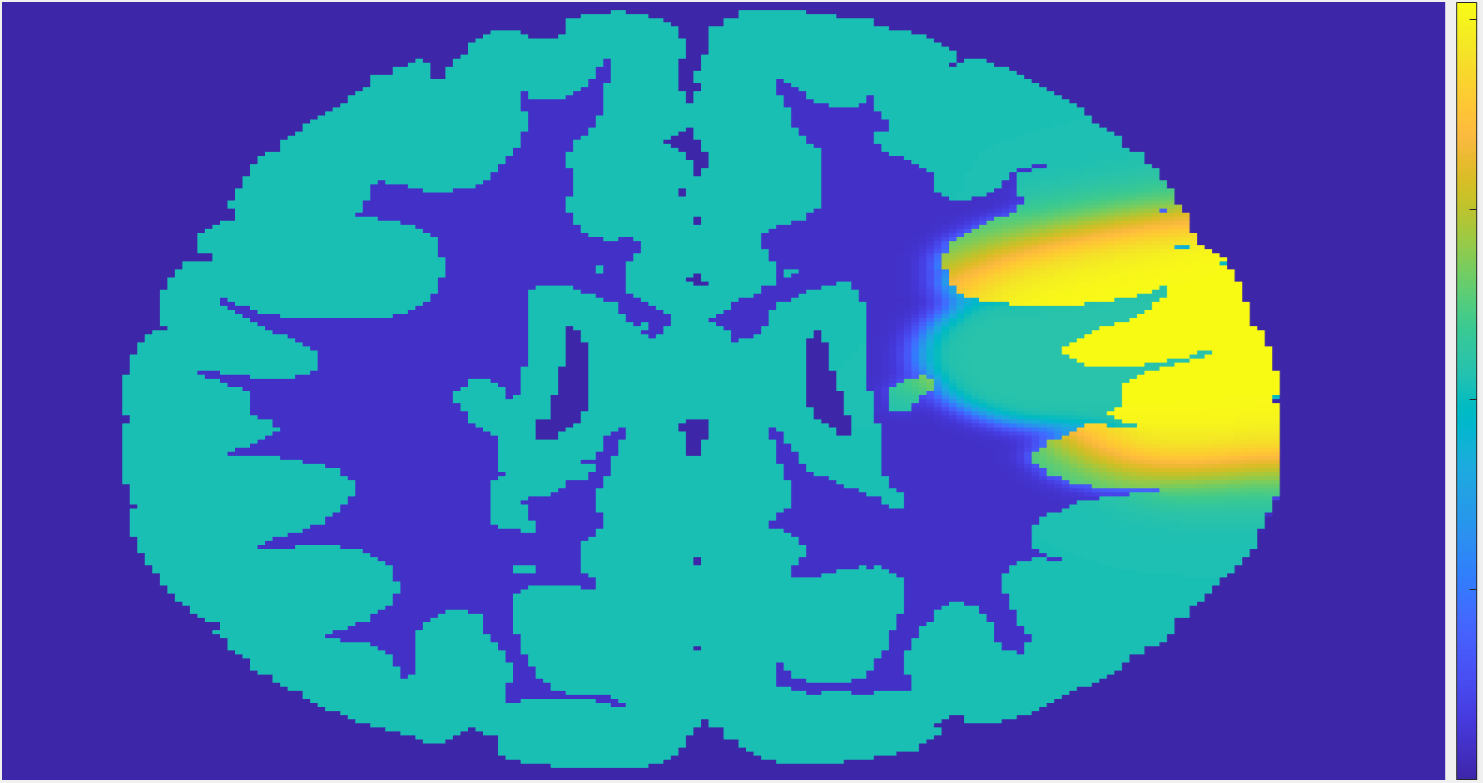}
         \caption{Snapshot at day 1222.}
         \label{fig:r} 
     \end{subfigure}
\\
     \begin{subfigure}[b]{0.49\textwidth}
         \centering
         \includegraphics[width=\textwidth]{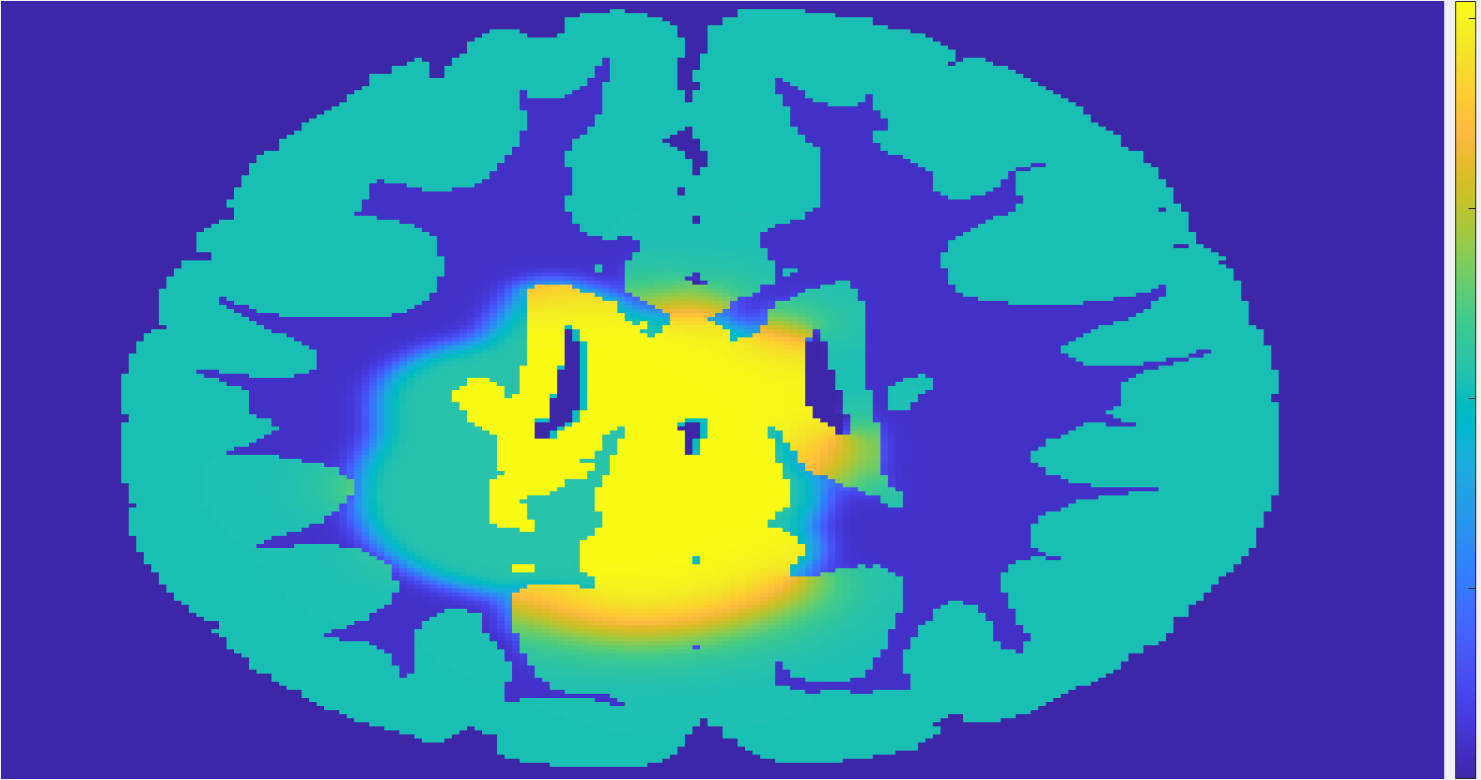}
         \caption{Snapshot at day 1356.}
         \label{fig:q}
     \end{subfigure}
     \hfill
     \begin{subfigure}[b]{0.49\textwidth}
         \centering         \includegraphics[width=\textwidth]{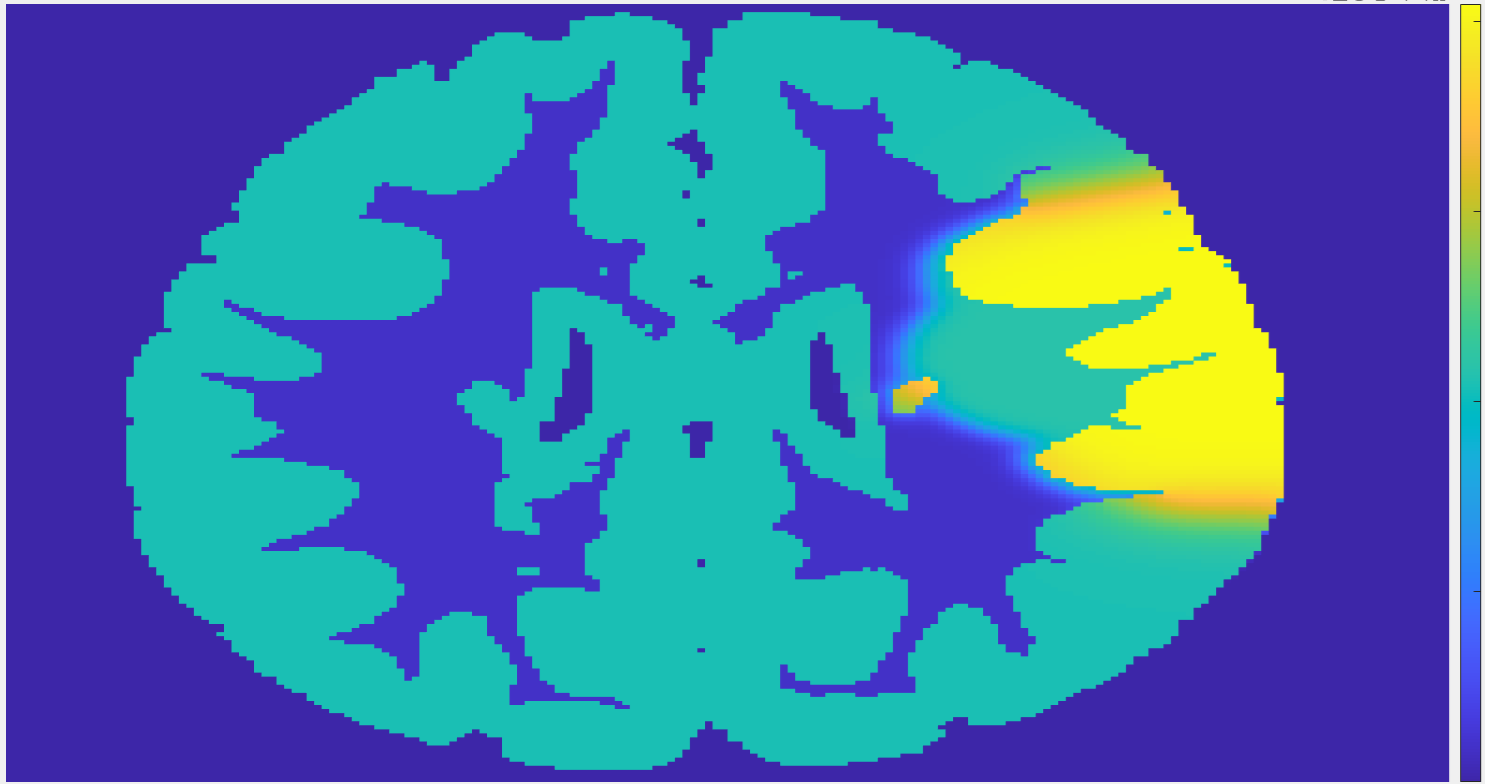}
         \caption{Snapshot at day 1356.}
         \label{fig:r} 
     \end{subfigure}\caption{\revision{Comparison of tumor growth starting in two different brain locations. Snapshots illustrate different growth rates in white and gray matter.}}
  \label{fig:comparison2}
\end{figure}

\revision{We present a comparison of tumor growth starting at two different locations in the brain. We employ the transverse view - a vertical cross-section of the brain. We start by denoting the white and gray matter on the cross-section (see Figure \ref{fig:matters}). The snapshots from the simulation are presented in Figures \ref{fig:comparison1}-\ref{fig:comparison2}. We can read from these figures that brain tumors grow faster in the white matter, and their final shape depends strongly on the initial location. }

{This single two-dimensional slice of the human head has around 30 percent of the gray matter and 70 percent of the white matter. According to our evaluation, we have theoretically predicted tumor velocity as
of \revision{$2\sqrt{\hat{D}\rho} =2\sqrt{(p_w D_w + p_g D_g) \rho} =2 \sqrt{(0.70\times 0.13 + 0.30\times 0.013) 0.025}= 0.048$
[mm days $^{-1}$]. 
This means that we have 0.48mm progression per 10 days.
The progression from day 150 (our initial state, where the tumor had a diameter of 1mm) to day 1356 (36-time steps of 33 days), presented in Figures \ref{fig:comparison1}-\ref{fig:comparison2}, shows the progression from 1mm to 62mm.
We identified the grid points with non-zero tumor cell density, and we found the largest distance between the two points as the tumor radius}.
This implies the tumor velocity of $\frac{62}{1356-150}=0.051$ [mm days $^{-1}$].}

\subsection{Three-dimensional glioblastoma brain tumor simulation}

We construct a 3D computational grid out of the slices of the MRI scan. We want to solve the Fisher-Kolmogorov diffusion-reaction equation with logistic growth \eqref{PDE} describing the brain tumor dynamics \cite{Logistic} over $\Omega\subset [0,200]^3$ [mm$^3$], and time interval {$I=(150,3500)$} [days]
and prescribed initial state
\begin{equation}\label{PDE3DBC}
{u(x_1,x_2,x_3;t_0)=0.1 \exp(-10((x_1-102)^2+(x_2-138)^2+(x_3-96)^2))}
\end{equation}
Following \cite{c25,Menze} we consider $\rho=2.5\times 0.01$ [days$^{-1}$]. $D_w=1.3\times 0.1$ [mm$^2$days$^{-1}$] \cite{c25}. $D_g=0.13\times 0.1$ [mm$^2$days$^{-1}$].

\begin{figure}[!htb]
     \begin{subfigure}[b]{0.49\textwidth}
         \centering
         \includegraphics[width=\textwidth,angle=180]{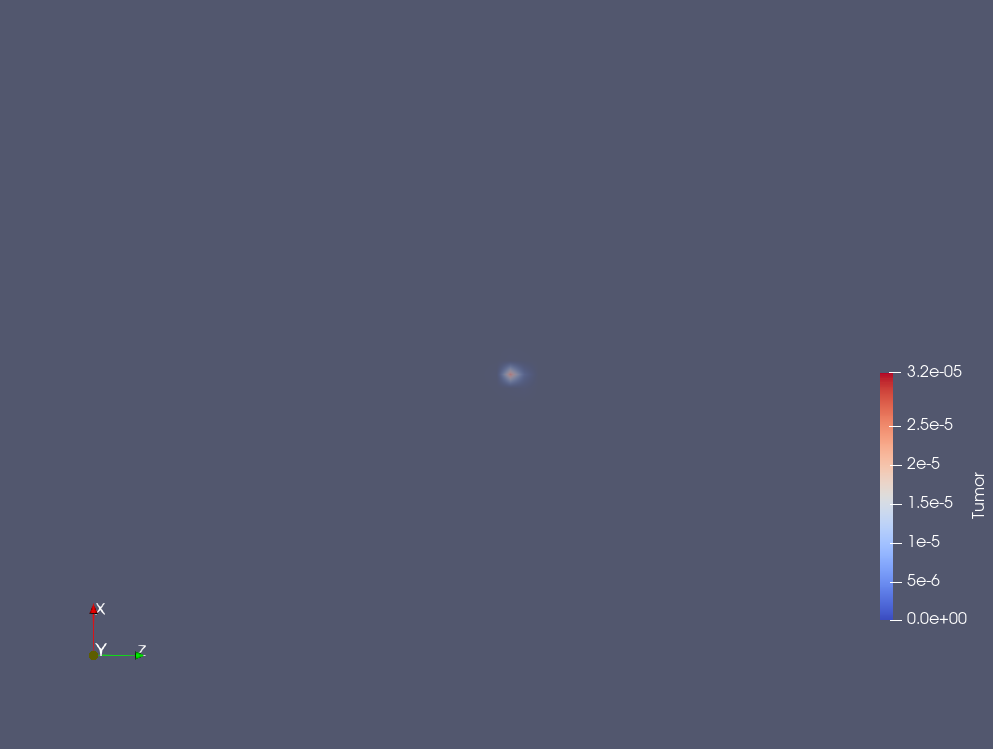}
         \caption{Sagittal view at day 150.}
         \label{fig:sa}
     \end{subfigure}
     \hfill
     \begin{subfigure}[b]{0.49\textwidth}
         \centering         \includegraphics[width=\textwidth,angle=180]{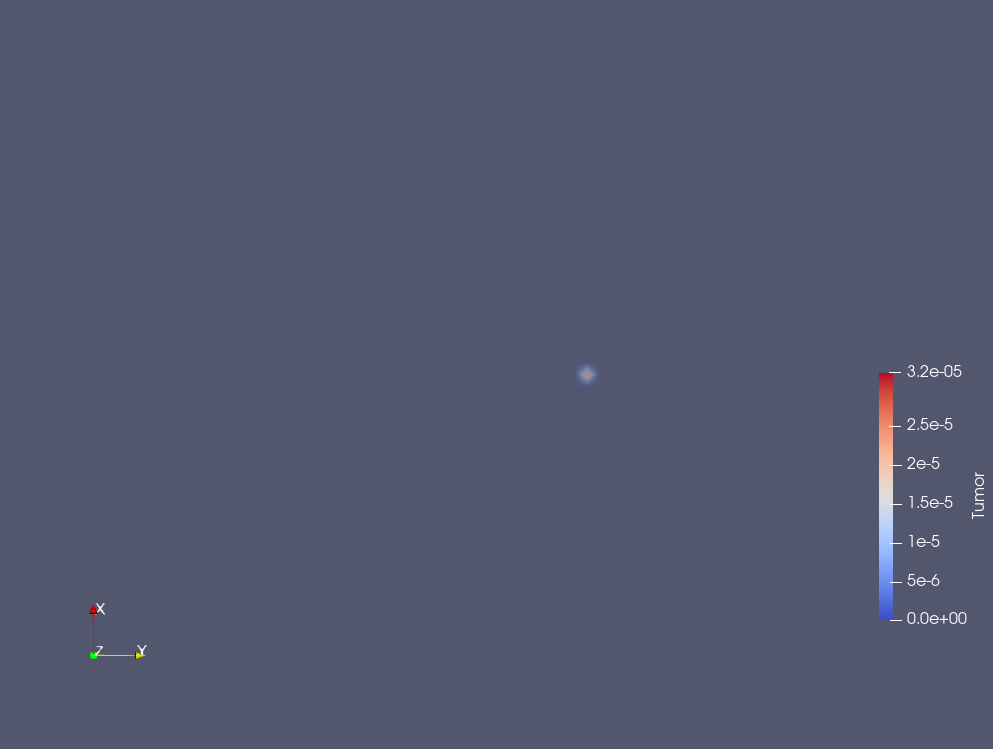}
         \caption{Coronal view at day 150.}
         \label{fig:ca}
     \end{subfigure}
     \begin{subfigure}[b]{0.49\textwidth}
         \centering
         \includegraphics[width=\textwidth,angle=180]{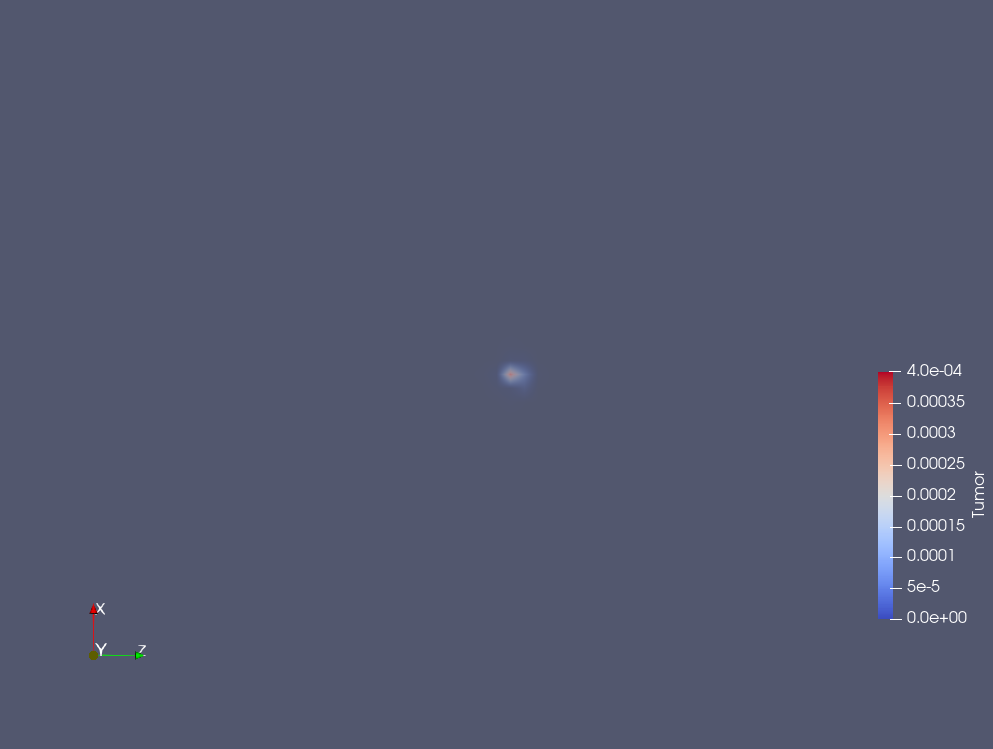}
         \caption{Sagittal view at day 315.}
         \label{fig:sb}
     \end{subfigure}
     \hfill
     \begin{subfigure}[b]{0.49\textwidth}
         \centering         \includegraphics[width=\textwidth,angle=180]{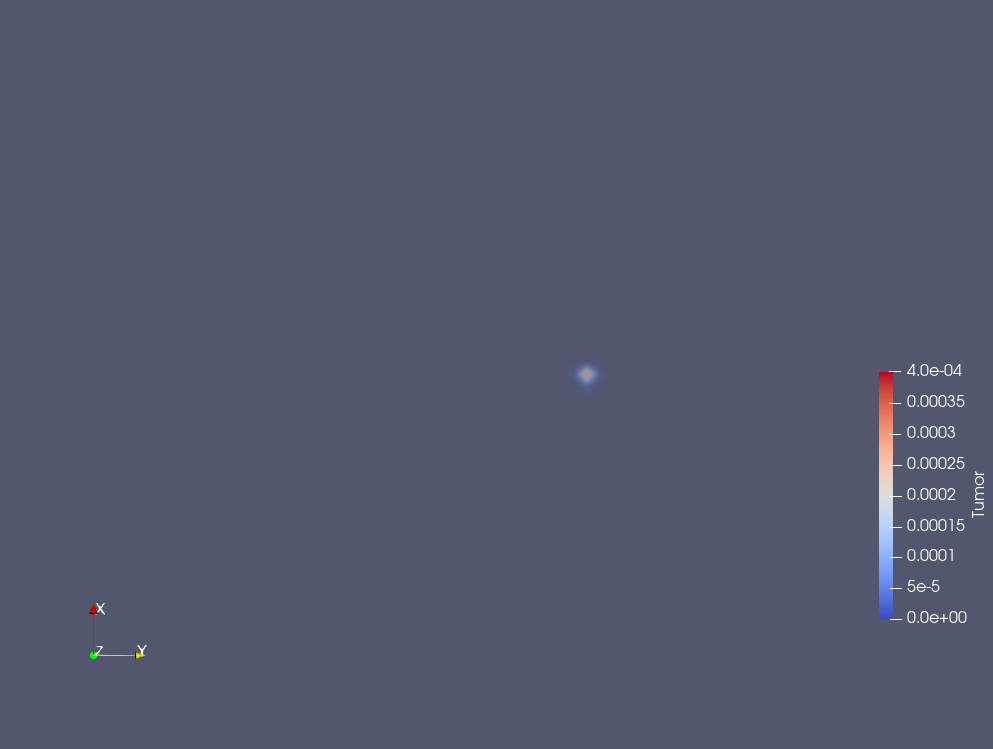}
         \caption{Coronal view at day 315.}
         \label{fig:cb}
     \end{subfigure}
     \begin{subfigure}[b]{0.49\textwidth}
         \centering
         \includegraphics[width=\textwidth,angle=180]{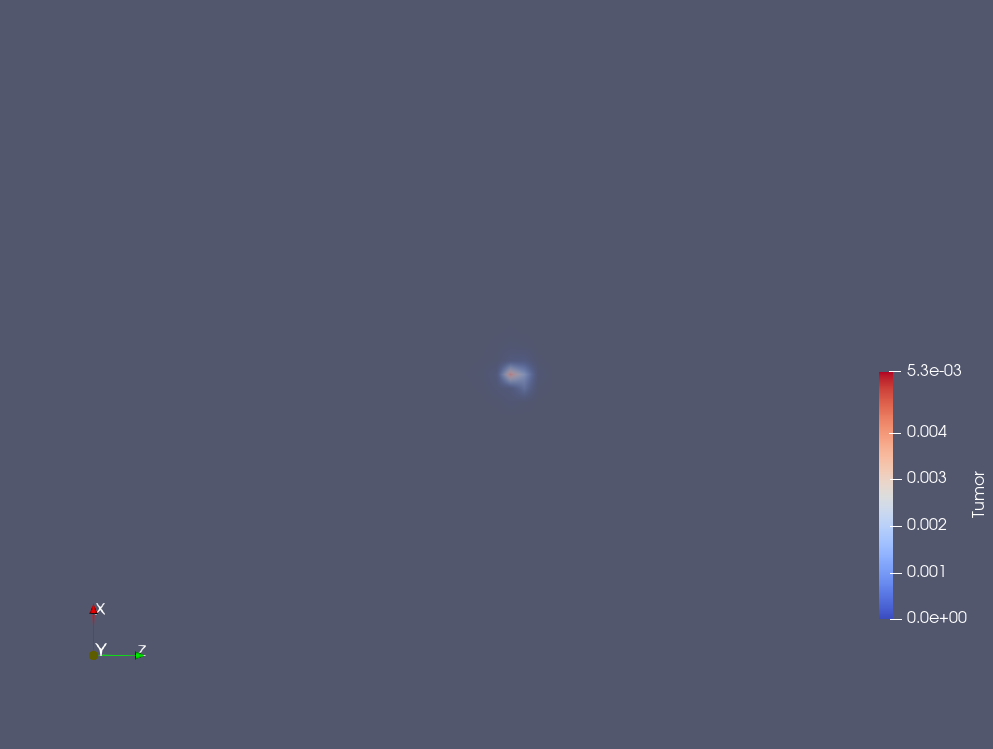}
         \caption{Sagittal view at day 480.}
         \label{fig:sa}
     \end{subfigure}
     \hfill
     \begin{subfigure}[b]{0.49\textwidth}
         \centering         \includegraphics[width=\textwidth,angle=180]{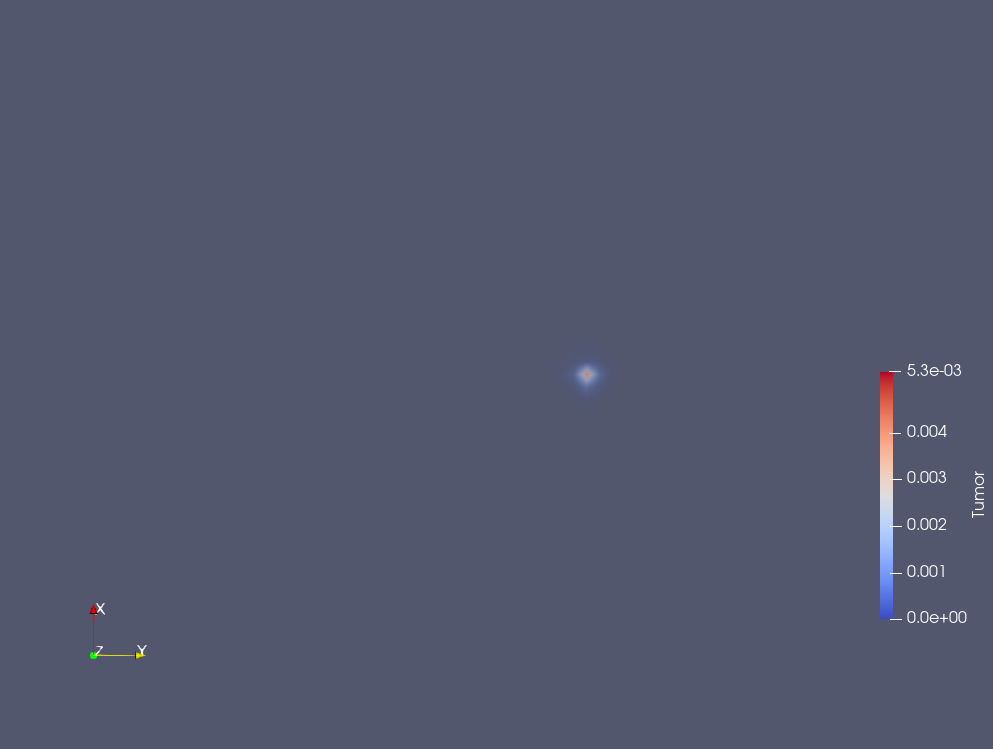}
         \caption{Coronal view at day 480.}
         \label{fig:ca}
     \end{subfigure}
         \caption{Saggital and coronal views of the brain tumor growing within white and gray matters on the cross-sections of the brain.}
         \label{fig:views}
\end{figure}

\begin{figure}[!htb]
     \begin{subfigure}[b]{0.49\textwidth}
         \centering
         \includegraphics[width=\textwidth,angle=180]{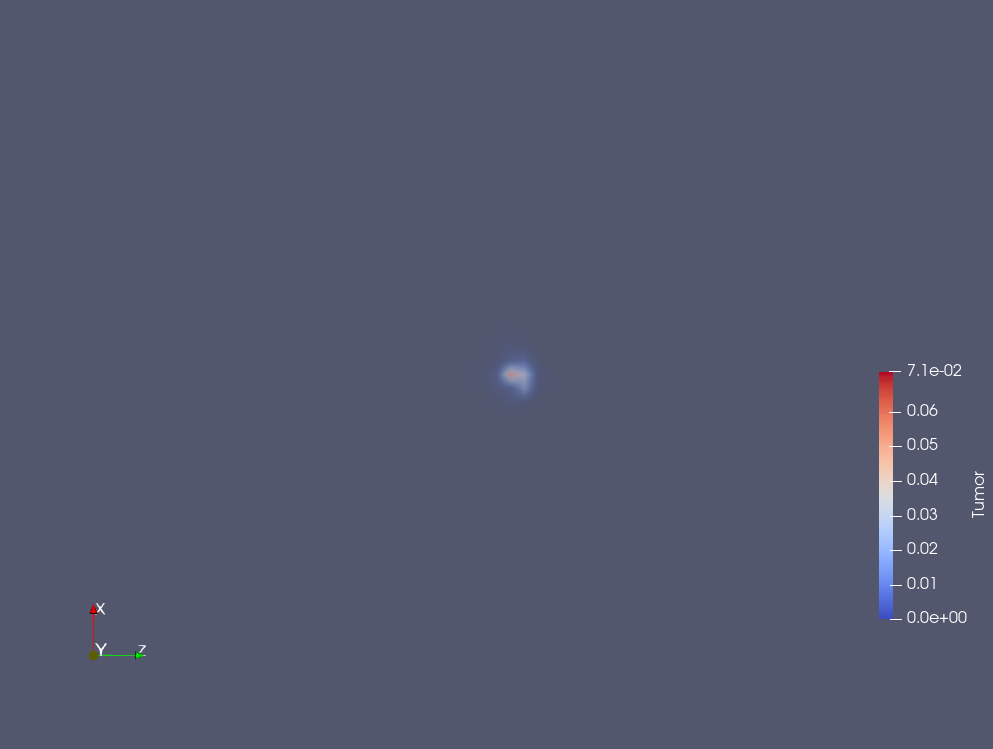}
         \caption{Sagittal view at day 645.}
         \label{fig:sa}
     \end{subfigure}
     \hfill
     \begin{subfigure}[b]{0.49\textwidth}
         \centering         \includegraphics[width=\textwidth,angle=180]{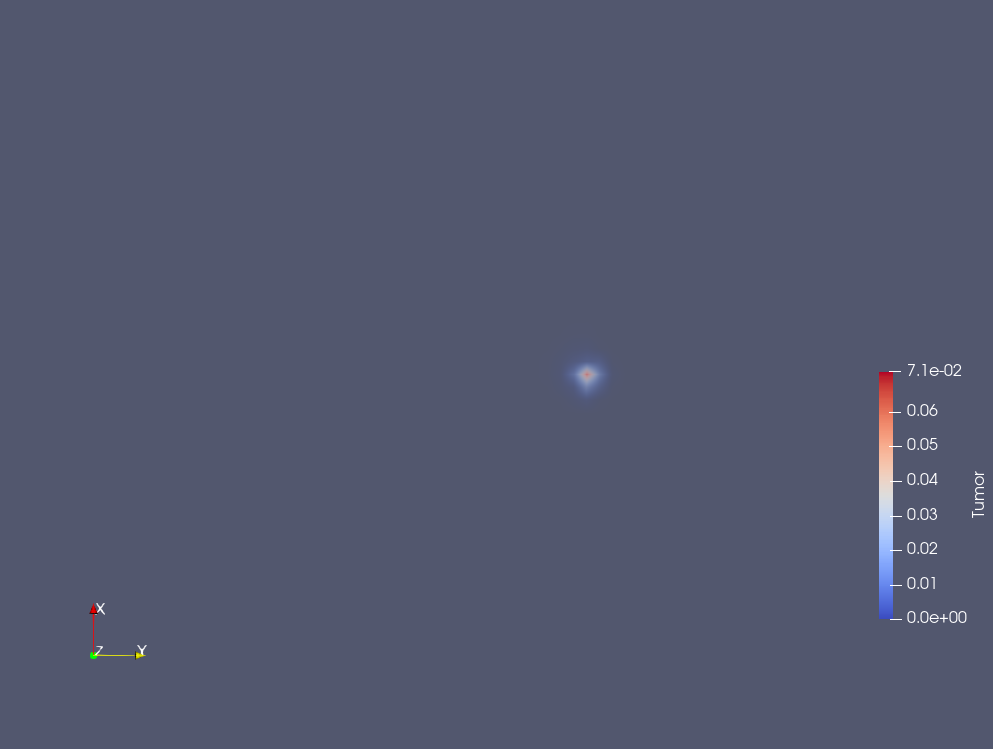}
         \caption{Coronal view at day 645.}
         \label{fig:ca}
     \end{subfigure}
     \begin{subfigure}[b]{0.49\textwidth}
         \centering
         \includegraphics[width=\textwidth,angle=180]{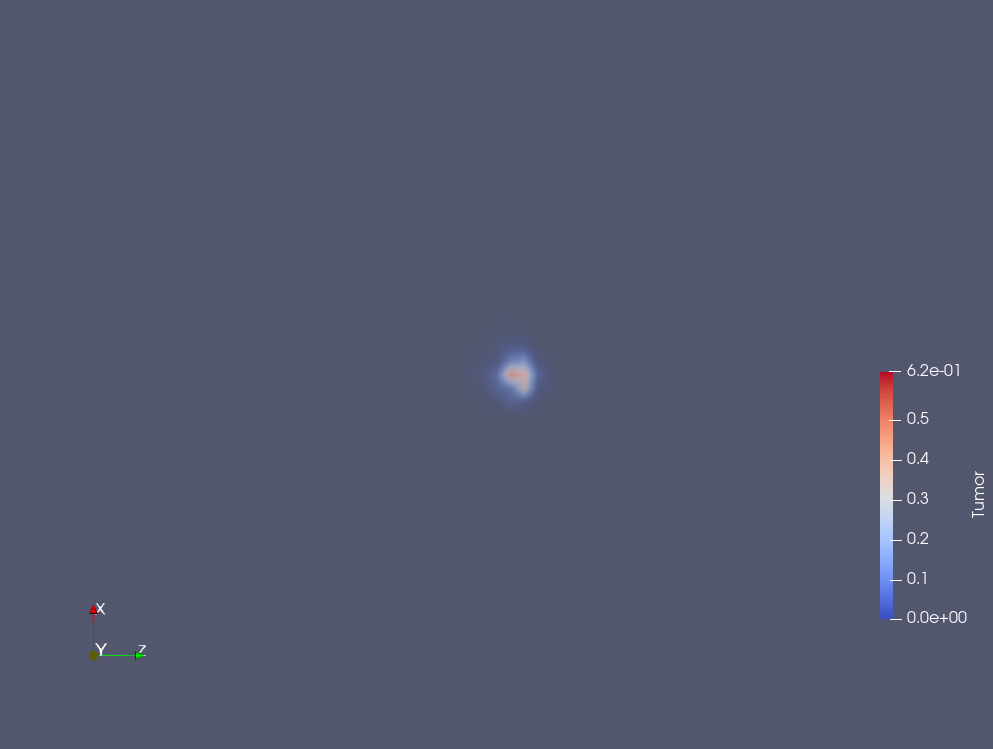}
         \caption{Sagittal view at day 810.}
         \label{fig:sa}
     \end{subfigure}
     \hfill
     \begin{subfigure}[b]{0.49\textwidth}
         \centering         \includegraphics[width=\textwidth,angle=180]{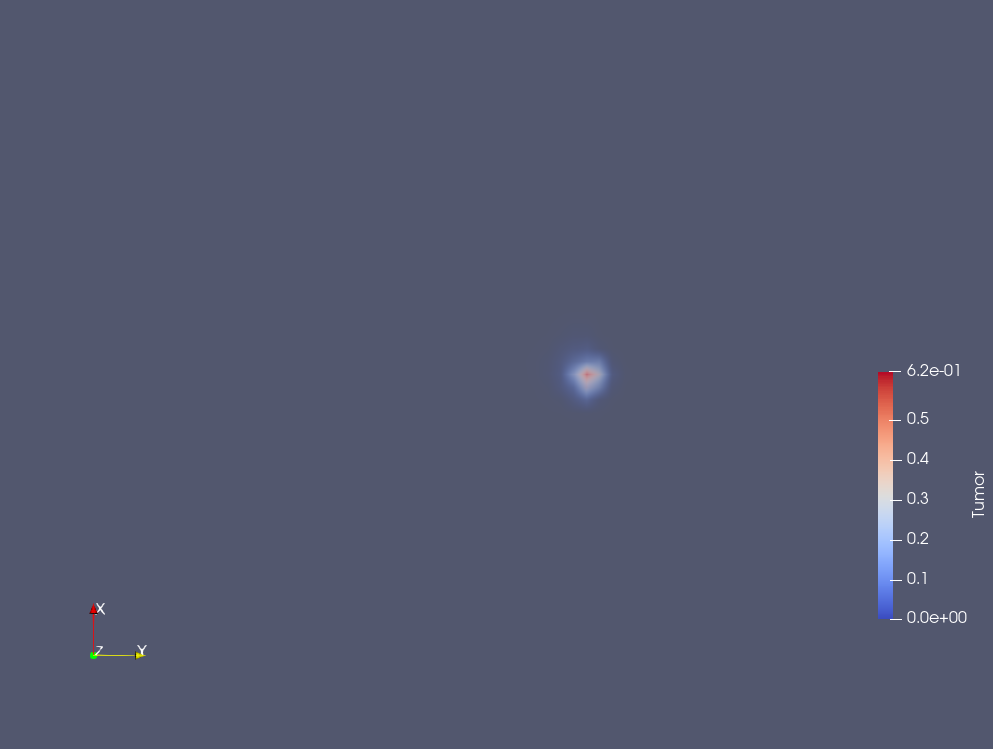}
         \caption{Coronal view at day 810.}
         \label{fig:ca}
     \end{subfigure}
     \begin{subfigure}[b]{0.49\textwidth}
         \centering
         \includegraphics[width=\textwidth,angle=180]{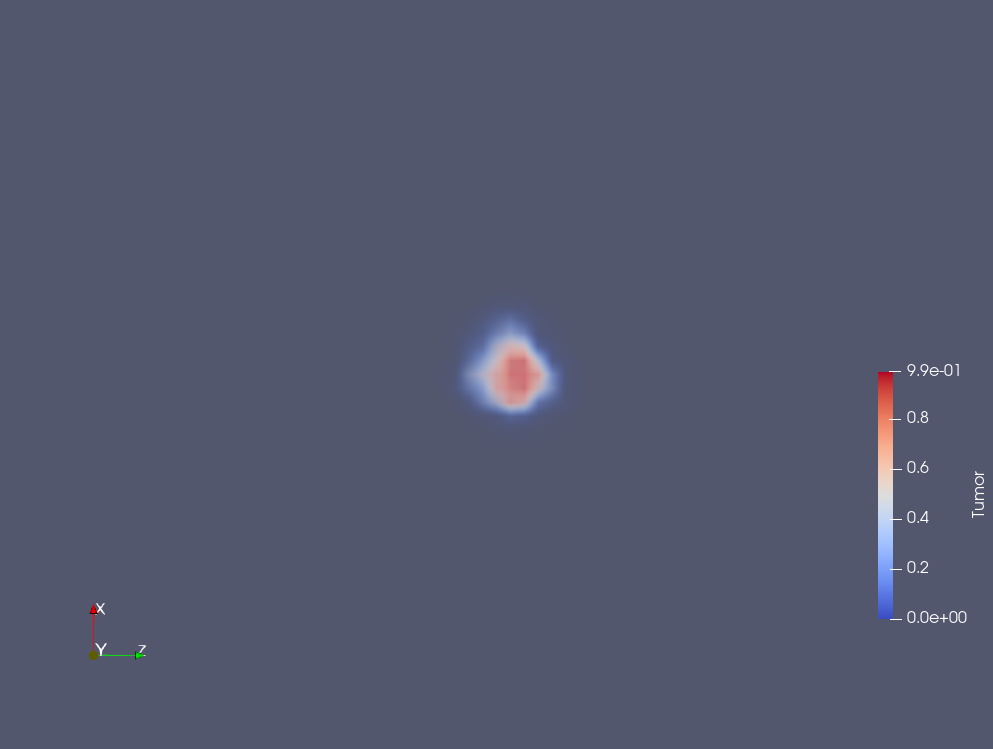}
         \caption{Sagittal view at day 975.}
         \label{fig:sb}
     \end{subfigure}
     \hfill
     \begin{subfigure}[b]{0.49\textwidth}
         \centering         \includegraphics[width=\textwidth,angle=180]{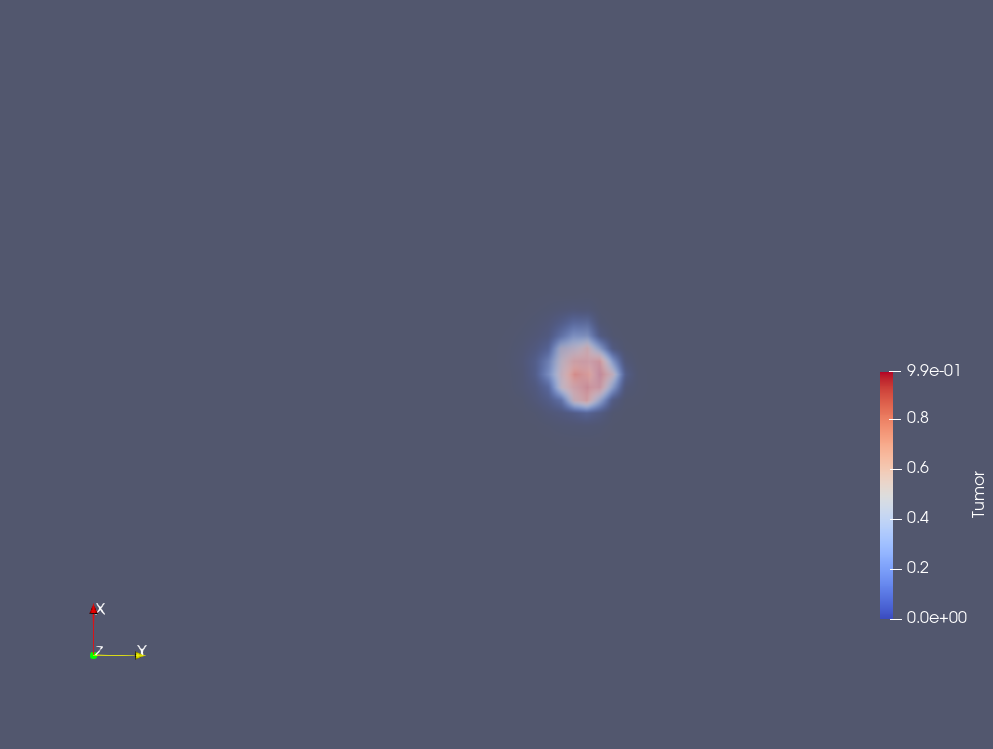}
         \caption{Coronal view at day 975.}
         \label{fig:cb}
     \end{subfigure}
         \caption{Saggital and coronal views of the brain tumor growing within white and gray matters on the cross-sections of the brain.}
         \label{fig:views2}
\end{figure}

\begin{figure}[!htb]
     \begin{subfigure}[b]{0.49\textwidth}
         \centering
         \includegraphics[width=\textwidth,angle=180]{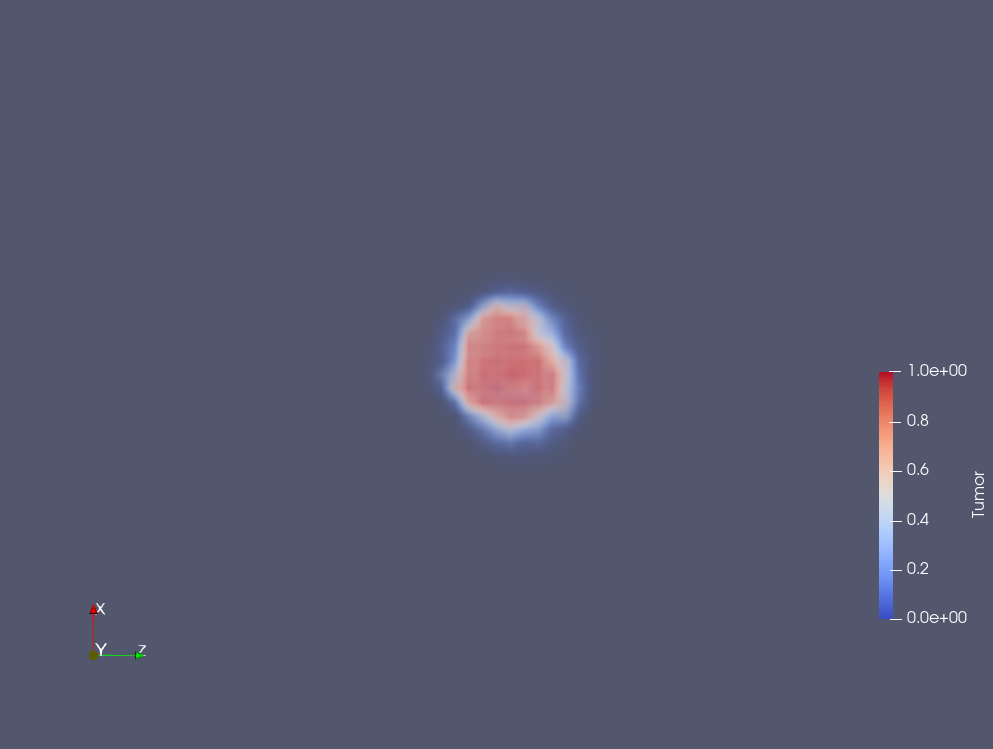}
         \caption{Sagittal view at day 1140.}
         \label{fig:sa}
     \end{subfigure}
     \hfill
     \begin{subfigure}[b]{0.49\textwidth}
         \centering         \includegraphics[width=\textwidth,angle=180]{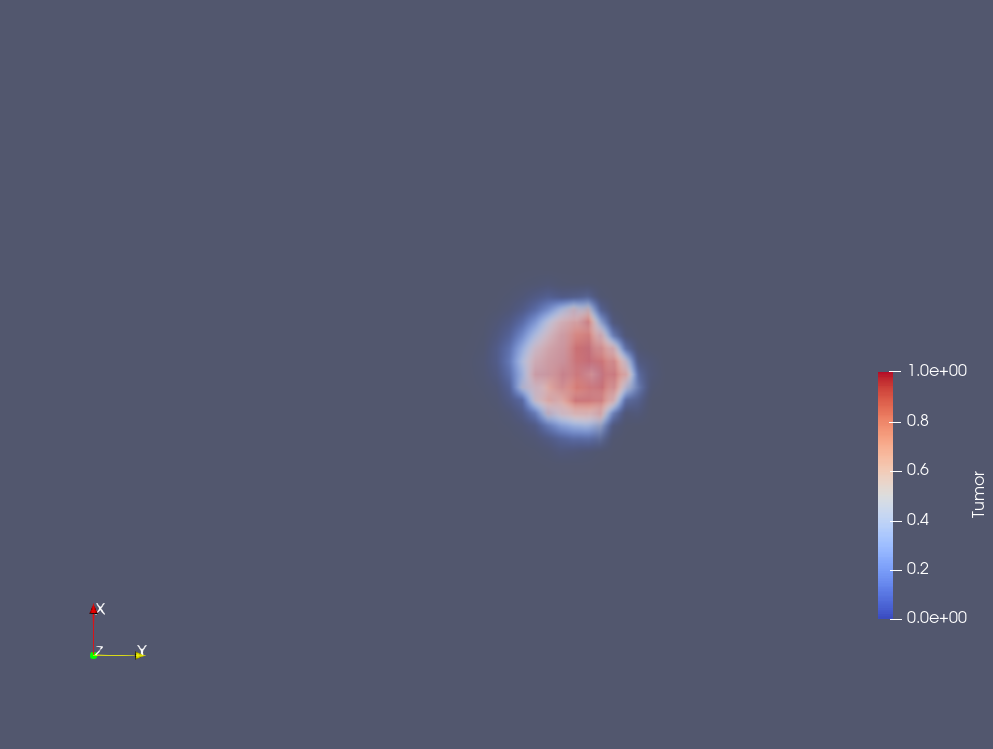}
         \caption{Coronal view at day 1140.}
         \label{fig:ca}
     \end{subfigure}
     \begin{subfigure}[b]{0.49\textwidth}
         \centering
         \includegraphics[width=\textwidth,angle=180]{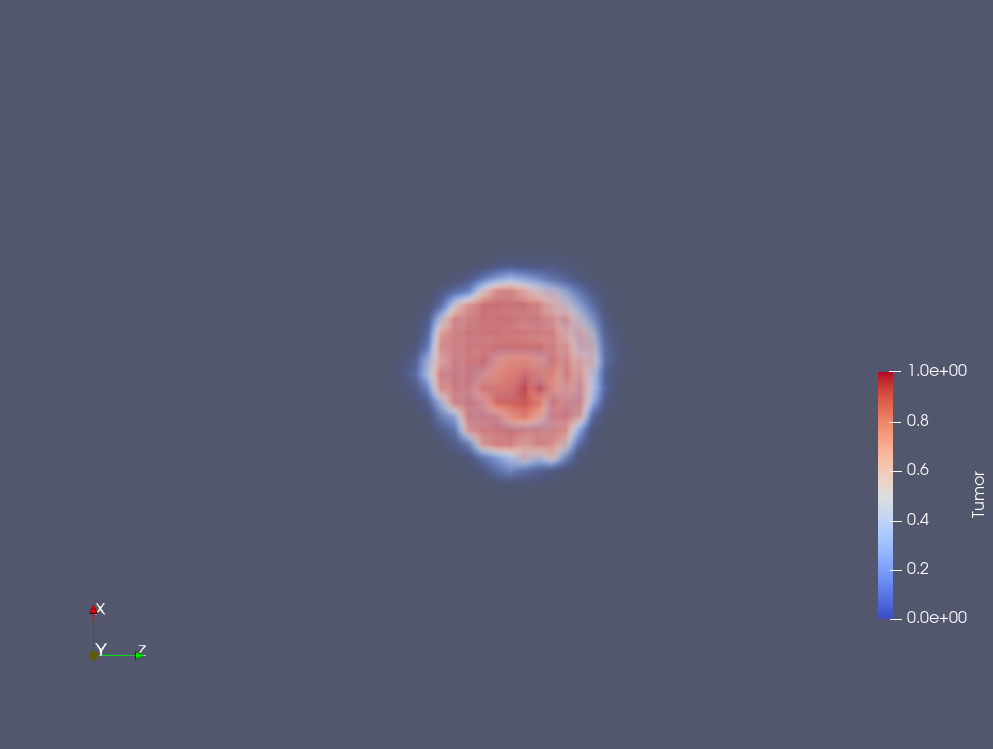}
         \caption{Sagittal view at day 1305.}
         \label{fig:sa}
     \end{subfigure}
     \hfill
     \begin{subfigure}[b]{0.49\textwidth}
         \centering         \includegraphics[width=\textwidth,angle=180]{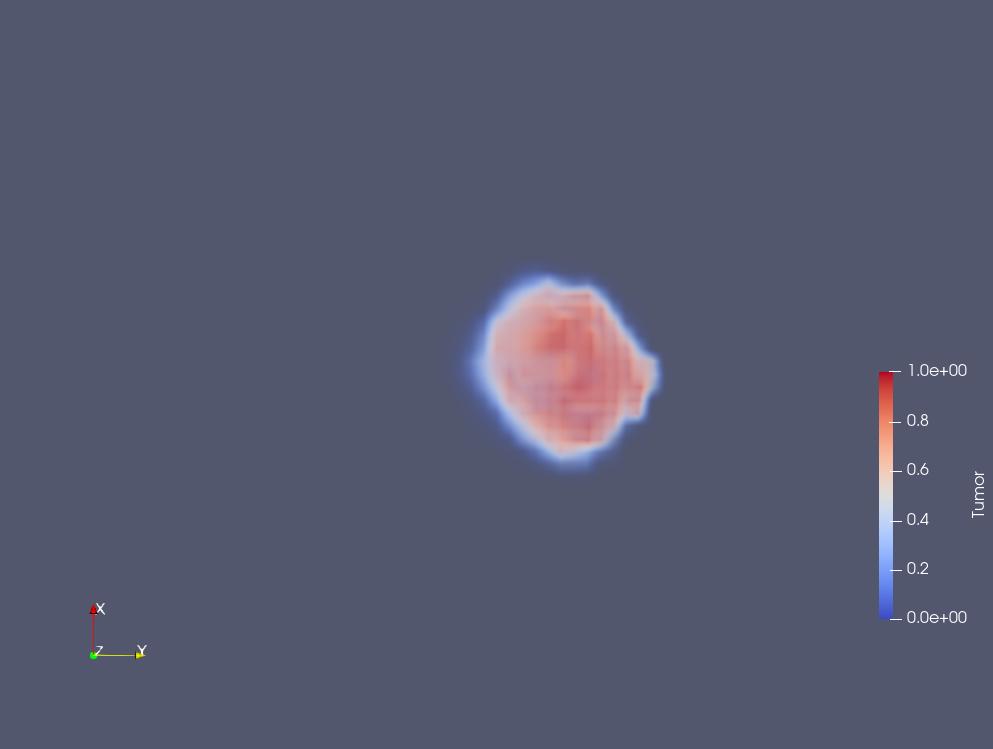}
         \caption{Coronal view at day 1305.}
         \label{fig:ca}
     \end{subfigure}
     \begin{subfigure}[b]{0.49\textwidth}
         \centering
         \includegraphics[width=\textwidth,angle=180]{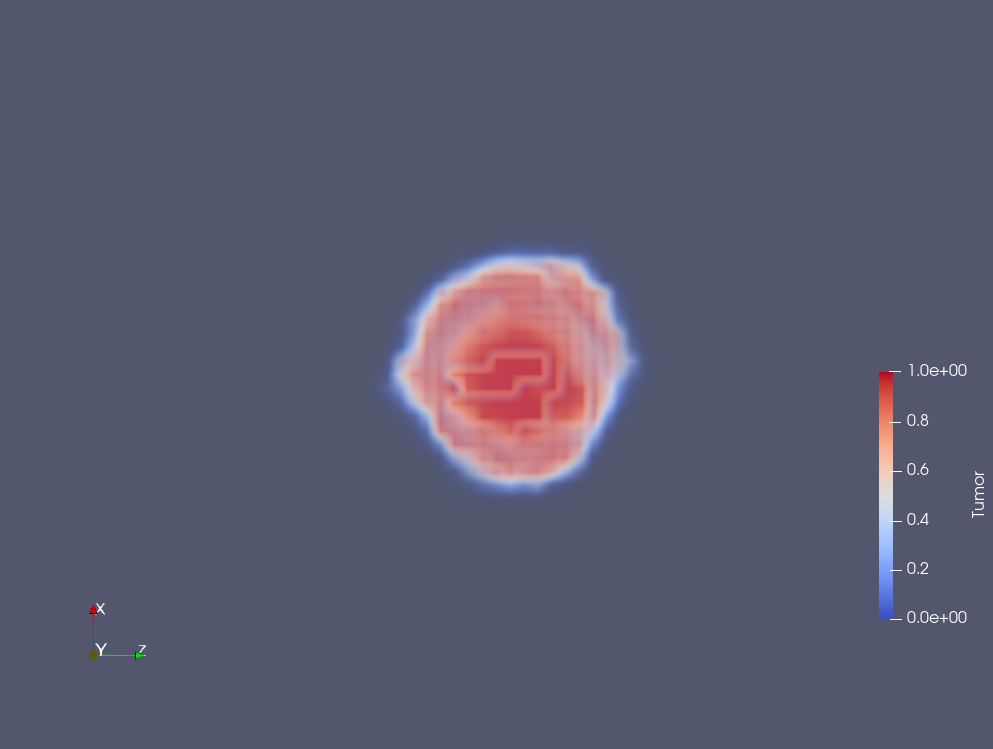}
         \caption{Sagittal view at day 1370.}
         \label{fig:sa}
     \end{subfigure}
     \hfill
     \begin{subfigure}[b]{0.49\textwidth}
         \centering         \includegraphics[width=\textwidth,angle=180]{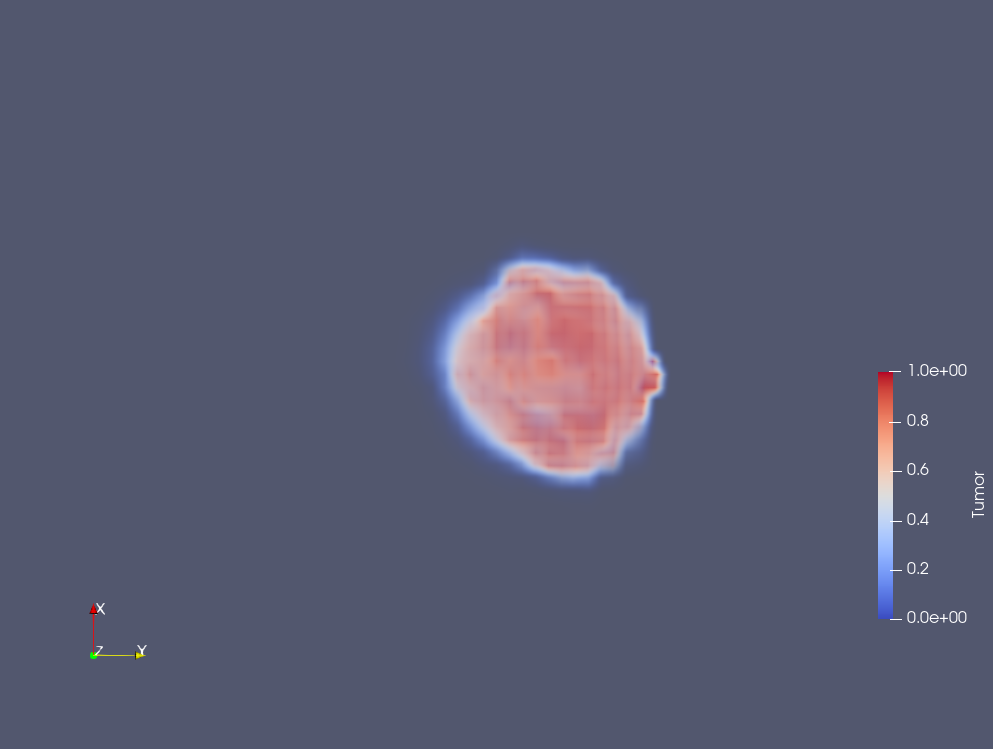}
         \caption{Coronal view at day 1370.}
         \label{fig:ca}
     \end{subfigure}
         \caption{Sagittal and coronal views of the brain tumor growing within white and gray matter on the cross-sections of the brain.}
         \label{fig:views3}
\end{figure}

We perform 100-time steps, starting from $t_0=150$ [days], with a time step $=33.5$ [days]. 
The spatial mesh consists of $128 \times 128 \times 128$ finite difference stencil points. 
The 100-time-steps exponential integrator method executed on a laptop with Win10, using MATLAB, with 11th Gen Intel(R) Core(TM) i5-11500H, 92 GHz and 32 GB of RAM, takes less than 10 minutes.
The selected snapshots from the simulation are presented in Figure \ref{fig:views}-\ref{fig:views3}. The visualization has been performed through the interface with ParaView \cite{ParaView}.

The percentage of the gray matter and white matter in the full three-dimensional MRI scan of the human head is around 85 percent of the gray matter and 15 percent of the white matter. This implies the following theoretical predicted tumor velocity as
of \revision{$2\sqrt{\hat{D}\rho} =$ 
$2\sqrt{(p_w D_w + p_g D_g) \rho} = 2 \sqrt{(0.15\times 0.13 + 0.85\times 0.013) 0.025}=0.055  $ 
[mm days $^{-1}$].
This means that we have 0.55mm progression per 10 days.}

\begin{figure}[!htb]
     \begin{subfigure}[b]{0.49\textwidth}
         \centering
         \includegraphics[width=\textwidth]{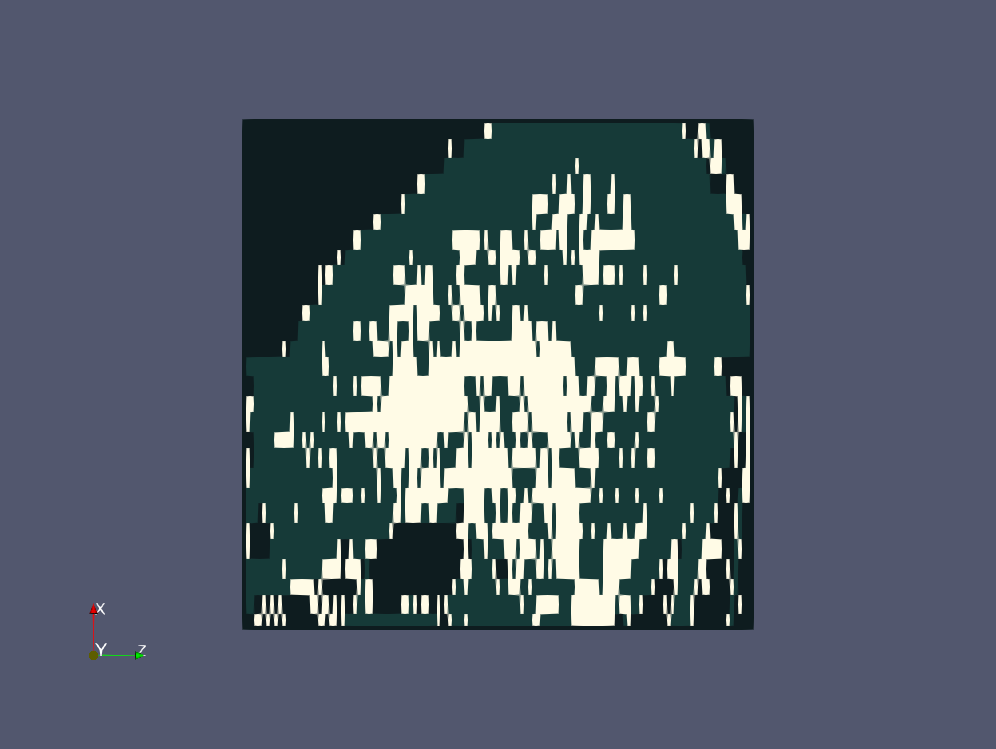}
         \caption{Sagittal view of the white and gray matter on the cross-section of the head.}
     \end{subfigure}
     \begin{subfigure}[b]{0.49\textwidth}
         \centering
         \includegraphics[width=\textwidth]{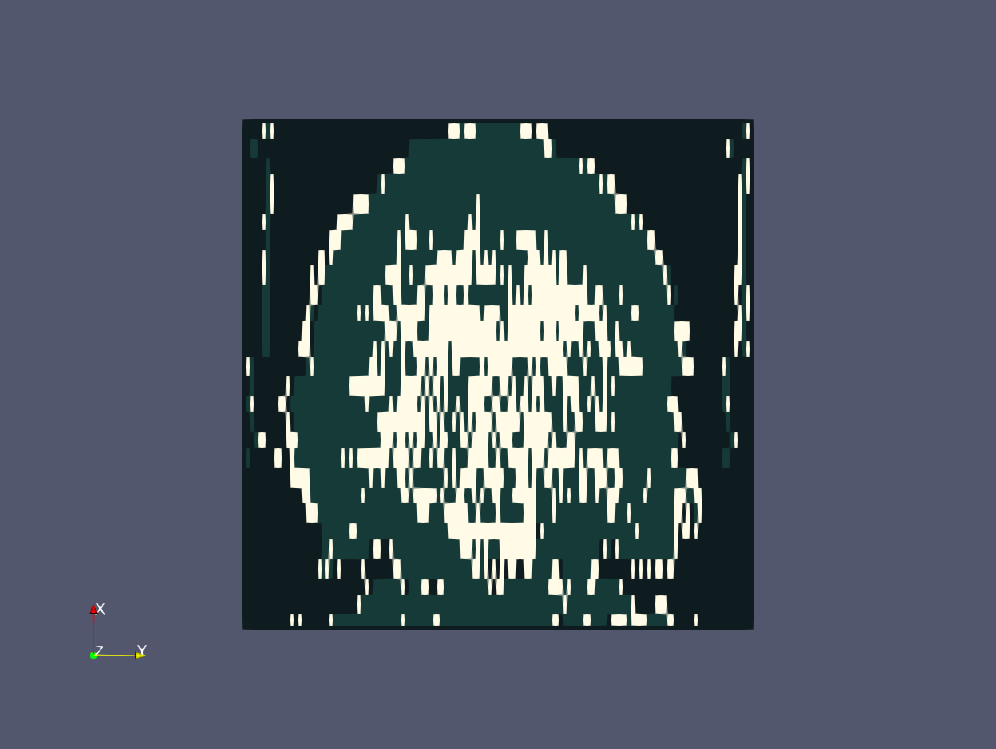}
         \caption{Coronal view of the white and gray matter  on the cross-section of the head.}
     \end{subfigure} \\

     \begin{subfigure}[b]{0.49\textwidth}
         \centering
         \includegraphics[width=\textwidth]{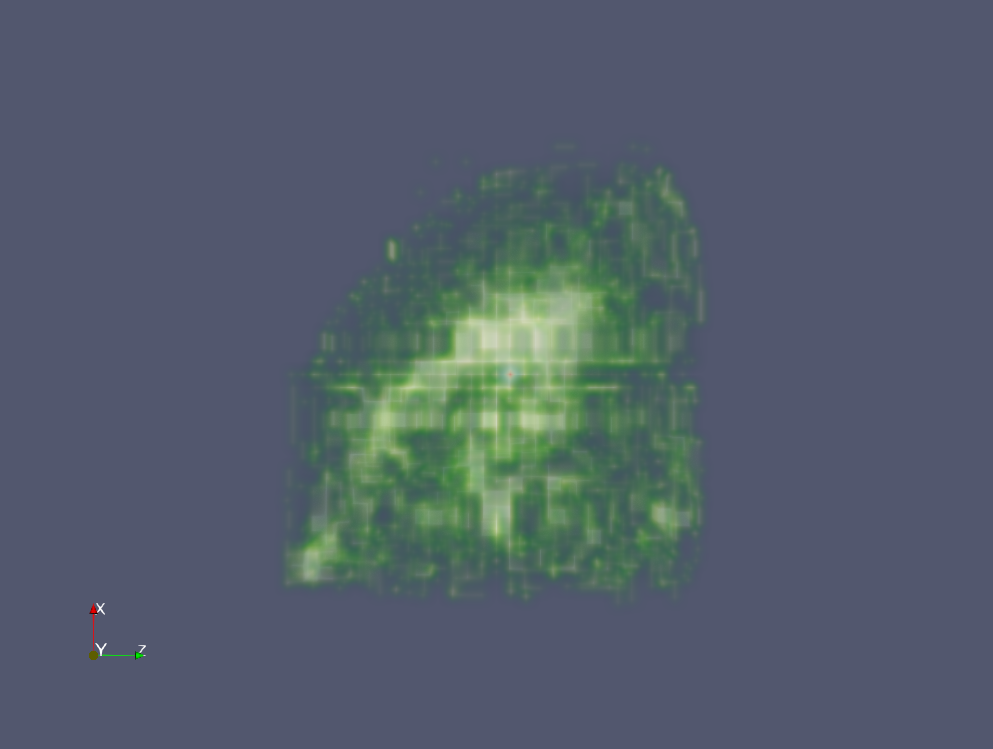}
         \caption{Sagittal view of the white and gray matter using Paraview visualization.}
     \end{subfigure}
     \begin{subfigure}[b]{0.49\textwidth}
         \centering
         \includegraphics[width=\textwidth]{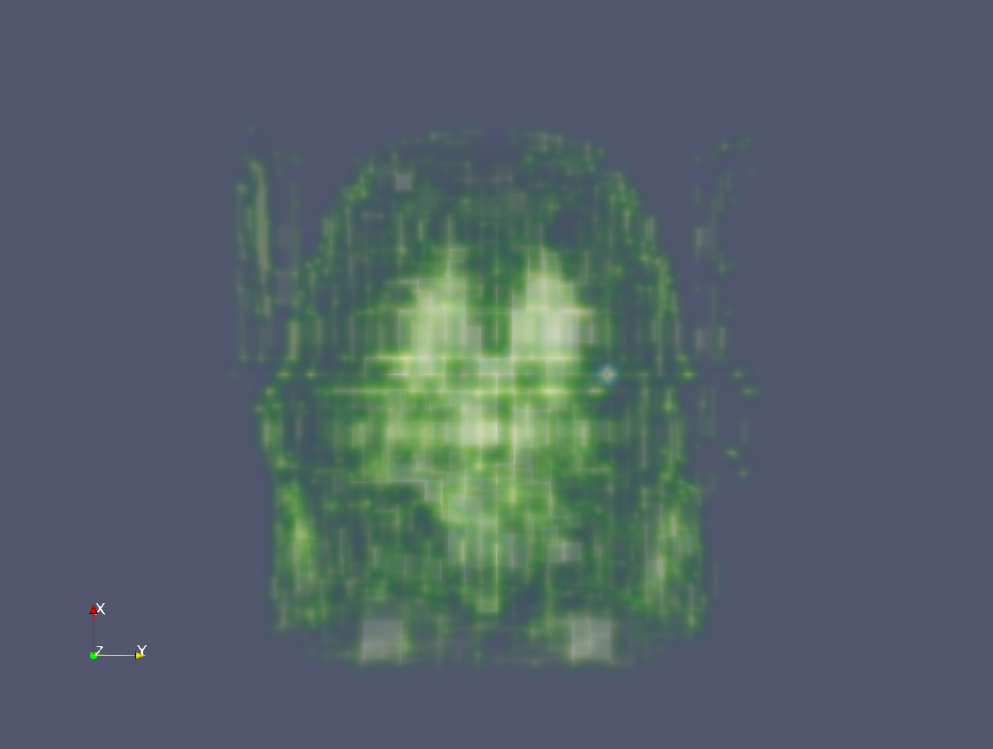}
         \caption{Coronal view of the white and gray matter using Paraview visualization.}
     \end{subfigure}
\caption{White and gray matters in the sagittal and coronal view at the cross-sections of the MRI scan data.}
         \label{fig:mattersx}
\end{figure}

\revision{The sagittal and coronal views of the white at gray matters at the cross-sections of the brain are presented in Figure \ref{fig:mattersx}. The sagittal and coronal views of the tumor growth simulation using Paraview visualization are presented in Figure \ref{fig:views}-\ref{fig:views3}. The progression from day 150 (our initial state, where the tumor had a diameter of 1mm) to day 975 shows the progression from 1mm to 45mm (measuring the tumor radius as a distance of two further points with non-zero tumor cell density) .
This implies the tumor velocity of $\frac{45}{975-150}=0.054$ [mm days $^{-1}$]. We have a good agreement between the experiments and the theory. }

\subsection{Comparison of the methods}
\revision{We first compare the execution times of 100 times steps of our brain tumor simulations calling the original routine \url{phiB} from \url{https://github.com/higham/expmv} library with our MATLAB implementation of the algorithm described in Section 2.2. Table \ref{tab:times} presents the execution times. Both MATLAB implementations are sequential; they use one core. We can read from this table that our
exponential integration routine is up to 3 times faster on $64^3$ and $96^3$ stencil points. For the $128^3$ stencil, the standard library cannot solve the problem (it "hangs out" at 3 percent of the computations for hours), and our library can solve $128^3$ stencil in less than 10 minutes.
We also compare our exponential integrators code with the Crank-Nicolson-based time discretization with the finite difference method.
We can read from this table that our code is up to two orders of magnitude faster. The Crank-Nicolson method cannot solve the problem on $96^3$ stencil on the laptop (it runs out of memory).}

\begin{table}[]
    \centering
    \begin{tabular}{|c|c|c|c|}
         \hline
      Stencil &  Exponential &  Exponential & Crank-Nicolson\\
      size &   generation &   integrators &  integrators  \\
& expmv.m & Sections 2.2-2.3 & \\
         \hline
$16^3$	&  	0,27 [s] & 0,38 [s] & 1,38 [s] \\
$32^3$	& 	1,96 [s] & 1,58 [s] & 42 [s] \\
$64^3$	& 	44 [s] & 14 [s] & 1807 [s] \\
$96^3$	& 	231 [s] & 96 [s] & - \\
$128^3$	& 	- & 507 [s] & - \\
         \hline
    \end{tabular}
    \caption{\revision{Execution times of EXPBrain code on the laptop using single core,  \url{https://github.com/higham/expmv} library, our algorithm from Sections 2.2-2.3, and Crank-Nicolson method.}}
    \label{tab:times}
\end{table}

\subsection{\revision{Convergence study}}

\begin{figure}[!htb]
  \centering
         \includegraphics[width=\textwidth]{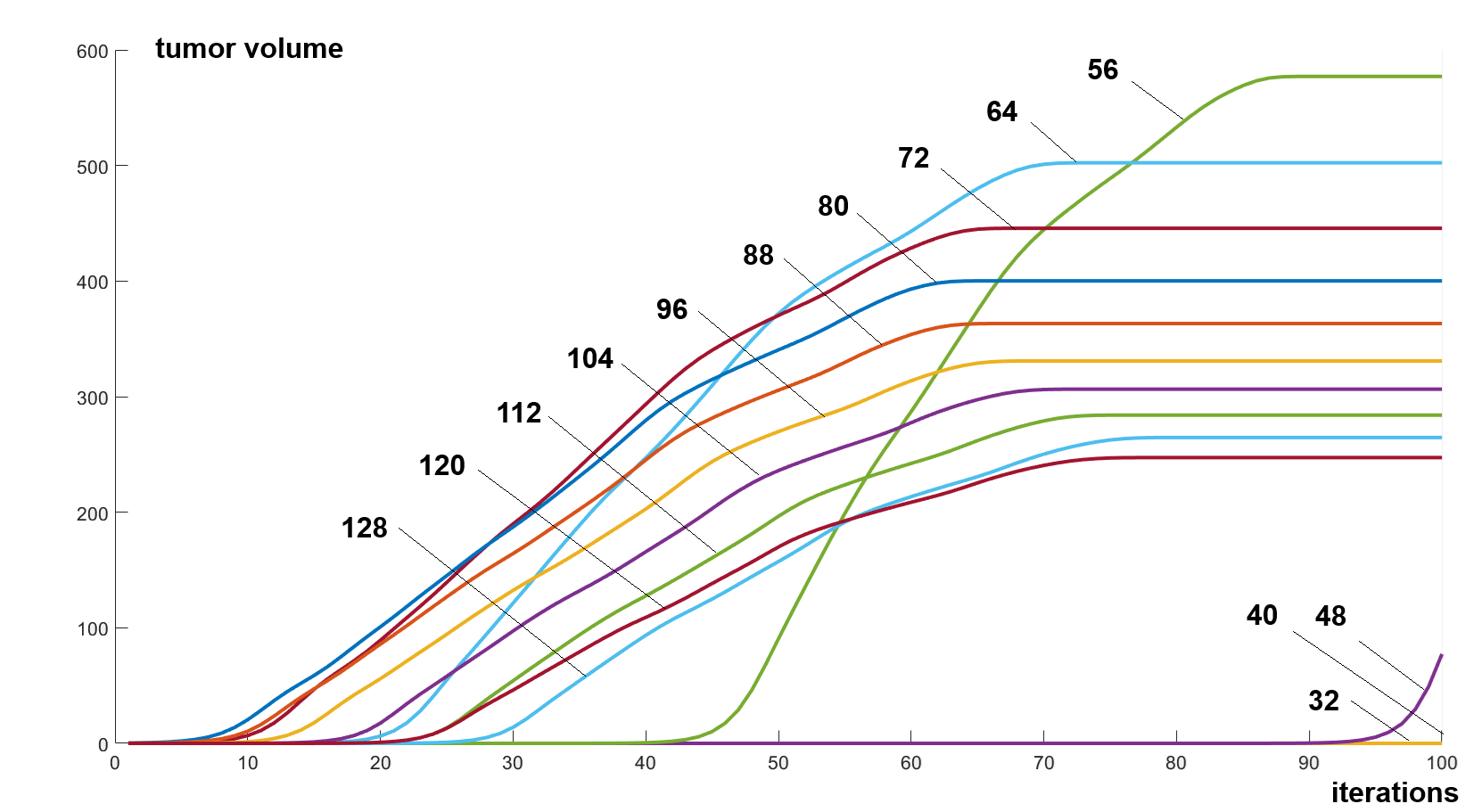}
         \caption{\revision{Convergence of the tumor growth simulations when varying the stencil size in each dimension. The first location of the tumor.}}
         \label{fig:convergence_all_scalled}
\end{figure}

\revision{We present the convergence study of the exponential integrator simulator, increasing the finite difference stencil in each direction. 
It is summarized in Figure \ref{fig:convergence_all_scalled}. 
The horizontal axis represents 100 iterations of the simulator. The vertical axis represents the tumor volume
computed as the summation of the tumor density at stencil points, divided by the number of stencil points. This estimation of the tumor volume is performed in every iteration.
Each line corresponds to a different dimension of the finite difference stencil.
In our simulations the tumor volume stops growing when the tumor fills the entire brain.
The dimensions of the MRI scan and the finite difference stencil are different. The MRI scan dimension remains constant and is usually larger than the finite difference stencil.
In our simulations we take values from the bitmap at the points of the finite difference stencil.
A lower resolution of finite difference stencil than the MRI scan results in some inaccuracies of the solution. Increasing the accuracy is illustrated in the convergence study in Figure \ref{fig:convergence_all_scalled}, measuring the $L^2$ norm of the tumor volume while increasing the finite difference stencil size.}

\begin{figure}[!htb]
  \centering
         \includegraphics[width=\textwidth]{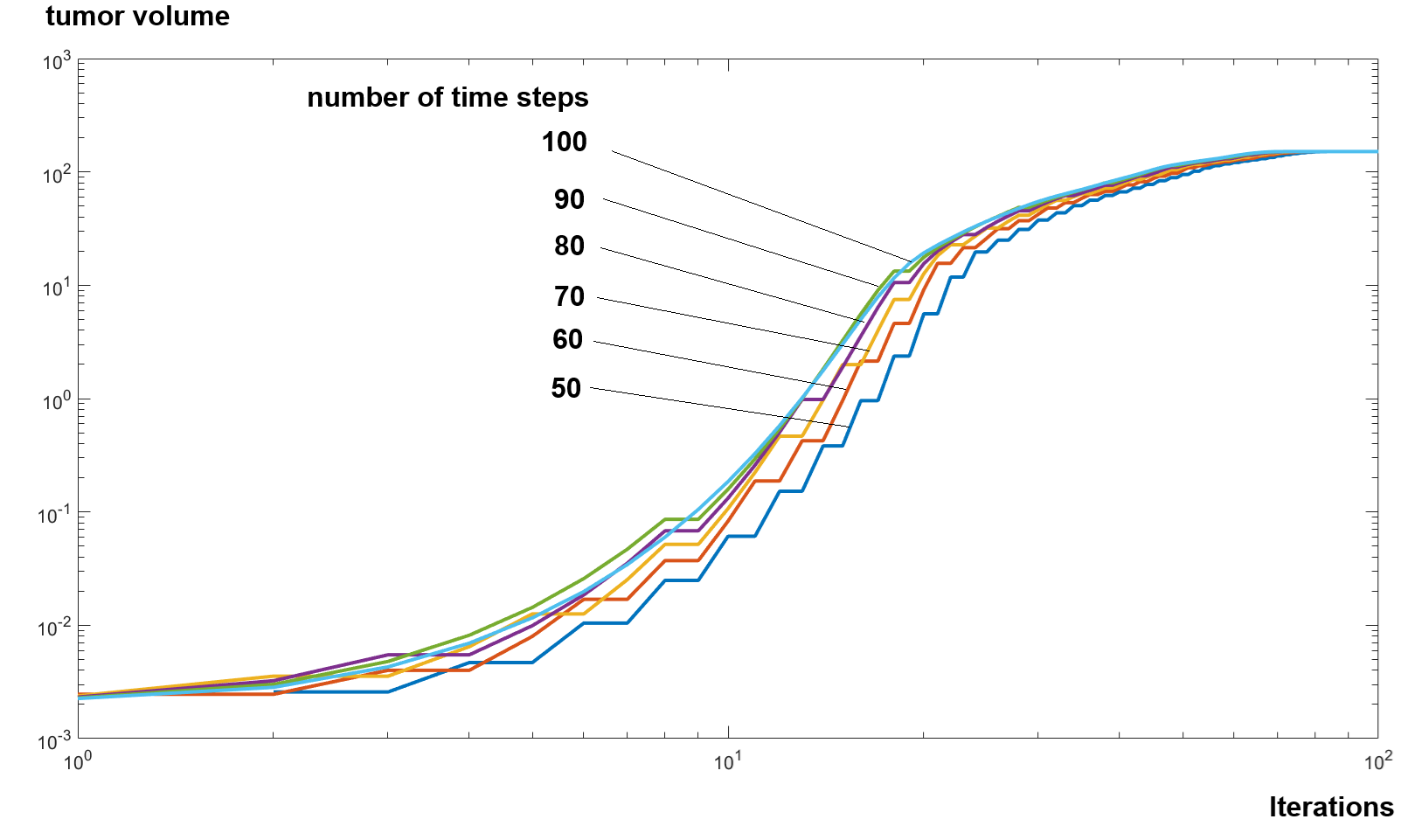}
         \caption{\revision{Convergence of the tumor growth simulations when varying the time step size.}}
         \label{fig:convergence_time}
\end{figure}

\revision{We also present the study of the convergence of the method while increasing the number of time steps executed in the time interval. We compare 50, 60, 70, 80, 90, and 100-time steps in Figure \ref{fig:convergence_time}. We see a good agreement of the results when varying the time step size. This time, we compute the tumor volume by summing up the tumor cells density at finite difference stencil points, but since we keep the stencil size fixed, we do not scale with respect to the number of stencil points.}

\clearpage
\section{Conclusions}
\revision{We presented a MATLAB code for predicting brain tumor growth on MRI scan data. It used the exponential integrators method in time and the finite difference method in space.}
Our two-dimensional simulator written in MATLAB can predict a one-year forward tumor evolution within 6 seconds using the resolution of 532$\times$565 pixels of a single two-dimensional slice of the human head.
\revision{We also developed a three-dimensional simulator with MATLAB. 
We proposed a new method for computing the exponential integrators and verified it on the brain tumor simulations.
The glioblastoma tumor growth model employed the Fisher-Kolmogorov diffusion-reaction model with logistic growth. }
The three-dimensional simulator on the laptop with Win10, using MATLAB, with 11th Gen Intel(R) Core(TM) i5-11500H, 2.92 GHz, and 32 GB of RAM, can perform a one-year tumor growth prediction.
\revision{The code, based on the provided MRI scan data of the human head and the initial tumor location, can predict the future evolution of the brain tumor. We can obtain the prediction for two years with an accuracy of $128^3$ on a laptop within less than 10 minutes.
Our method is faster than state-of-the-art exponential integrators library \cite{al2011computing}, it can process larger stencils, and it is two orders of magnitude faster than the Crank-Nicolson method with finite difference discretization, which actually can simulate a maximum size of the stencil $64^3$ within half an hour.}
For both two-dimensional and three-dimensional simulations, we set the parameters of the simulator based on the average percentage of white and gray matter.
With good agreement, we have compared the average velocity of the tumor growth resulting from the numerical experiments with the velocity predicted by the theory.
In future work, we plan to analyze the concurrency of the exponential integrators algorithm and develop a parallel version of the simulator that targets GPGPU computing cards. Using the parallel simulator, we will run large simulations on high-resolution human brain data. Another future direction of our research is the incorporation of the drug delivery term into our model \cite{McDaniel2013}. In particular, we plan to focus on the chemotherapy modeling. We also plan to address the data assimilation problem using the supermodeling technique \cite{PASZYNSKI2022214,SIWIK2022114308}.
\revision{The future work will also include incorporating adaptive mesh refinements \cite{hp3d}. We will try to incorporate our graph-grammar approach \cite{GG1,GG2,GG3} into the computational method of exponential integrators.}

\section*{Acknowledgements}
This work was partially supported by the program ``Excellence initiative - research university" for the AGH University of Krak\'ow. Judit Mu\~noz-Matute has received funding from the European Union's Horizon 2020 research and innovation programme under the Marie Sklodowska-Curie individual fellowship grant agreement No. 101017984 (GEODPG).
\revision{We would like to thank the anonymous reviewers for their valuable suggestions.}

\section*{\revision{Credit statements}}
\revision{{\bf Magdalena Pabisz}: Conceptualization, Investigation, Methodology, Software, Validation, Visualization, Writing – original draft, {\bf Judit Muñoz-Matute}: Conceptualization, Formal analysis, Funding acquisition, Investigation, Methodology, Software, Validation, Writing – review \& editing, {\bf Maciej Paszy\'nski}: Conceptualization, Data curation, Formal analysis, Funding acquisition, Project administration, Supervision, Writing – original draft, Writing – review \& editing. }


\end{document}